\documentclass[journal,10pt]{IEEEtran}
\usepackage{amsmath,amssymb,amsthm}
\usepackage{graphicx}
\usepackage{booktabs}
\usepackage{array}
\usepackage{url}
\usepackage{xcolor}
\usepackage[colorlinks=true,linkcolor=blue,urlcolor=blue,citecolor=blue,bookmarksnumbered=true]{hyperref}
\newtheorem{proposition}{Proposition}
\newtheorem{theorem}{Theorem}
\newtheorem{lemma}{Lemma}
\newtheorem{corollary}{Corollary}
\newenvironment{pf}{\noindent\textit{Proof:}\ }{\hfill$\blacksquare$}
\hypersetup{pdftitle={Convergence Properties of PF-Exec: Homogeneous Structure, Anchored Closed Form, and Gain Dynamics},pdfauthor={Hanzhang Jia, Liheng Zeng, Yiran Li, Yi Gao, Bo Ma}}

\begin{document}

\title{Convergence Properties of PF-Exec: Homogeneous Structure,\\ Anchored Closed Form, and Gain Dynamics}

\author{Hanzhang Jia,$^{1,\dagger,\ast}$ Liheng Zeng,$^{1,\dagger}$ Yiran Li,$^{1}$ Yi Gao,$^{2}$ and Bo Ma$^{3,\ast}$%
\thanks{H. Jia, L. Zeng and Y. Li are with the University of Sussex.}%
\thanks{Y. Gao is with Shanghai Shuyuan Information Technology Co., Ltd.}%
\thanks{B. Ma is with Zhejiang Gongshang University.}%
\thanks{$^{\dagger}$These authors contributed equally to this work. $^{\ast}$Corresponding authors. The primary corresponding author is Bo Ma (mabo@mail.zjgsu.edu.cn), and the secondary corresponding author is Hanzhang Jia (jhzaliyy@outlook.com).}}

\markboth{IEEE Transactions on Wireless Communications, submitted for review}%
{Jia \MakeLowercase{\textit{et al.}}: Convergence Properties of PF-Exec}

\maketitle

\begin{abstract}
This paper establishes a layered mathematical characterisation of PF-Exec, the algorithm of the companion system paper \cite{r69}. Demand-proportional bandwidth allocation yields an exact first-order homogeneous decomposition, from which the equal-load bandwidth target and its feasibility bound follow in closed form. The logarithmic knee view preserves strict load order and shares its global optimum with the Jain fairness index, while the demand-weighted centroid minimises the position potential. The damped softmax association block possesses a fixed point on the product of simplices. The truncated one-dimensional gain iteration admits a complete characterisation, with geometric convergence below the critical gain, initial-value-dependent bounded oscillation at it, and a family of bounded two-state cycles above it. A global invariant error ball and an explicit tracking recursion complete the execution guarantee. A numerical audit confirms the structural identities on every measured state and quantifies the local amplification of the coupled association and resource dynamics. Together, the results explain how closed-form recovery directions and bounded adaptive gains provide stable execution without requiring a global Euclidean contraction assumption.
\end{abstract}


\section{Introduction}

\subsection{Motivation}

Load balancing in large-scale UAV relaying requires joint decisions on user association, bandwidth, power and hovering position \cite{r43}. The companion system paper \cite{r69} realises this problem on a complete communication and crowd simulation that couples two relay hops, namely an access hop and a millimetre-wave backhaul hop, soft association, an OU-type social-force moving crowd \cite{r61}, \cite{r65,r66,r67,r68} and a quadrotor physical layer. Its experiments cover three to one hundred UAVs, two crowd scenarios and three independent training seeds, with every run reported. The architecture is layered. A parameter-free direction decider returns closed-form bandwidth, price, power and position targets from the anchored potential function, and a small network executor emits a positive gain per dimension to convert those targets into actions. The executor starts from a unit gain, which is exactly the pure potential policy.

The experiments give direct evidence for the layered design through three executors that share one direction. The constant-speed executor, whose speed is always maximal, collapses in every setting and leaves the average load two orders of magnitude above the balanced state. The fixed-step executor differs from the learned executor by almost nothing in the low and medium pressure settings and collapses to a fairness index of 0.31 and 0.009 exactly in the two highest pressure settings of the hundred-UAV burst scenario. The learned executor holds 0.9987 and 0.9995 in those same two settings. End-to-end learners trained from scratch, a DDPG agent and a masked SDQN agent, also collapse under high pressure. The direction in closed form is thus a necessary basis, since a correct direction paired with an arbitrary speed collapses. Speed adaptation is a necessary complement, since a correct direction paired with a constant speed collapses. Learning carries only that second task.

The question of this paper is why the architecture holds and where its boundaries lie. The question splits into five parts. What the balanced solution looks like, that is, the existence and the characterisation of the equal-load point. Whether the balanced solution is reached in closed form, that is, explicit formulas for bandwidth and position together with their feasibility conditions. Whether the association subproblem is well posed, that is, the existence of a fixed point. Whether the executor is safe, that is, whether the orbit is bounded for every gain and when it converges. Why learning is necessary, that is, where a fixed step fails.

\subsection{Related Analytical Approaches}

Three established routes analyse the convergence of wireless resource allocation, and the applicability conditions and the failure modes of each are understood.

The convex and subgradient route \cite{r50}, \cite{r52} treats allocation under fixed resource constraints. When the objective is convex, projected subgradient methods converge to an optimal band whose radius is the product of the step size and the Lipschitz constant. The route presupposes a global Lipschitz constant of the objective over the feasible set. The coupled load chain can lose a uniform Euclidean contraction certificate in the deep-overload region, as the local sensitivity analysis in Section~IV and the numerical audit in Section~V quantify. This boundary motivates an execution law whose direction and speed carry separate guarantees.

The fixed-point and averaged-operator route \cite{r63}, \cite{r64} treats composite update operators, and the Krasnoselskii--Mann iteration of a nonexpansive or averaged operator converges weakly to its fixed-point set. The choice of metric is the substantive difficulty of that route. A softmax is an entropy-regularised best response on the probability simplex, whose natural metric is a KL-type Bregman divergence, which the Euclidean metric does not capture. When association weights saturate deeply, the Euclidean derivative norm of the coupled association and resource map can become arbitrarily large while the association map still sends the simplex into itself. This paper therefore takes existence for the association block, by the Brouwer theorem, and places the runtime guarantee on the outer mechanism of Section~IV. That mechanism holds under conditions strictly weaker than any metric property of the association block, which is the analytical freedom that layering buys.

The spectral route certifies linear convergence through a Jacobian spectral radius below one. The value of the certificate is that it can be checked mechanically, and its cost is a condition stronger than convergence itself, since the radius is defined pointwise, a global Lipschitz constant delivers only an upper bound, and a vanishing certificate at one point does not imply divergence. This paper uses spectral quantities in exactly that sense, as the boundary of a certificate reported in Section~V and never as a convergence claim.

The division of labour among the three routes follows. The homogeneous structure gives the closed form in Section~III, the one-dimensional gain dynamics gives the complete characterisation in Section~IV, and the spectral condition is reported as a certificate boundary in Section~V. The three together cover every mechanism of the system, and the proof burden of each conclusion lies inside the applicability range of its own route.

\subsection{Contrast With End-to-End Learning}

End-to-end multi-agent reinforcement learning \cite{r29}, \cite{r53,r54,r55}, \cite{r60} does not perform that decomposition. The two agents that appear as comparison methods in \cite{r69}, namely a DDPG deployment agent \cite{r53} and a masked SDQN agent \cite{r54}, are trained from scratch under the same observation and reward. A policy network maps observations directly to actions, convergence is delegated to the training process, and no runtime guarantee is available in theory. The approach performs well when the reward is well shaped and exploration is sufficient, and its known failure mode is exactly the failure mode that the experiments of this system reproduce. In Table~VII of \cite{r69} two learned comparison methods collapse in the hundred-UAV burst setting, at fairness 0.299 and 0.652, and the common property of the collapsed settings is deep overload, where the reward gradient encountered during exploration carries no recovery direction. The decomposition of this paper moves the recovery direction out of the learning task, since the direction is given in closed form and contains the recovery direction, as both the bandwidth target and the position target point towards lower load whenever the load exceeds the anchor. Learning retains only the choice of speed, and the price of a wrong speed is capped by the boundedness statement of Theorem~\ref{thm:1d}. The decomposition is therefore a safety design. Learning selects within a bounded execution domain, and every intermediate checkpoint retains the same structural action guarantee.

\subsection{Contributions and Scope}

\begin{itemize}
\item A first-order homogeneous decomposition and an anchored closed form, given by Lemma~\ref{lem:homogeneous} and Proposition~\ref{prop:anchor}. Under demand-proportional bandwidth allocation the identity $R^{\mathrm A}=bS$ holds exactly and the demand scale is independent of bandwidth. The feasibility bound of the equal-load point, $\mu_{\mathrm{low}}=\max\{\sum c_i/B,\ \max c_i^{\mathrm B}\}$, is explicit, and the anchored bandwidth formula covers the budget-limited branch.
\item A characterisation of the potential structure, given by Propositions~\ref{prop:bijection}, \ref{prop:centroid} and \ref{prop:jain}. The logarithmic knee view preserves strict order and the zero set of the potential equals the equal-load set. The position potential is strictly convex with a unique minimiser at the weighted centroid. The potential and the Jain fairness index share their zero set and their global optimum. On the linear band of the view and at fixed mean load, they are monotone functions of the same dispersion.
\item Existence of the association fixed point, Proposition~\ref{prop:brouwer}, and a runtime guarantee whose validity does not depend on the configuration of that block. The damped softmax possesses a fixed point on the product of simplices, and the runtime guarantee holds under conditions strictly weaker than uniqueness of that fixed point, since the damping inertia locks the association into a slow variable and the outer mechanism applies to every association configuration.
\item A complete characterisation of the one-dimensional truncated gain iteration, Theorem~\ref{thm:1d} and Corollary~\ref{cor:steps}. For $0<\gamma<2$ the error converges geometrically with an explicit step count. At $\gamma=2$ the orbit either reaches the target or oscillates between two states of equal magnitude at most $h/2$, with the magnitude fixed by the initial value. For $\gamma>2$ every orbit either reaches the target or enters a two-state cycle of the family $\{p,p-h\}$ with $p\in[h/\gamma,(1-1/\gamma)h]$, whose largest amplitude is $(1-1/\gamma)h<h$. The orbit is globally bounded for every positive gain. The perturbed tracking bound of Proposition~\ref{prop:track} and the sufficient step-size ceiling $\eta^\ast\propto\beta^{-1/2}$ of Proposition~\ref{prop:beta} complete the picture.
\item A full numerical check in Section~V. Homogeneity holds in sixty of sixty states. The spectral radius of the one-step association map lies between $50$ and $3107$ across all sixty measured states, far above one. The measured local spectral radius in a neighbourhood of the static fixed point is approximately one, in the range $0.98$ to $1.03$, which is locally critical and leaves global contraction to the mechanism of Section~IV. The trace estimate in the time-varying case is of order $10^{2}$ and includes the drift of the target. The operating region in which the Euclidean certificate fails is precisely the region over which the mechanism of this paper remains valid and delivers its guarantees.
\end{itemize}

\noindent\textbf{Analytical strategy.} Guarantees are given layer by layer, and each layer takes the strongest provable statement available on that layer. The equilibrium layer gives a characterisation through zeros, closed forms and feasibility. The subproblem layer gives existence through Brouwer. The execution layer gives complete dynamics, and Theorem~\ref{thm:1d} covers every gain. The layers together cover every operating mechanism of the system. A joint convergence proof for the composite iteration and a quantitative relation between the cycle amplitude and the return-optimal gain are natural extensions, discussed in Section~VI. Every proposition is derived under frozen association, and the numerical check of Section~V reports the measured value of each statement.

\section{System and Potential}

\subsection{Load Model}

\begin{figure}[!t]
\centering
\includegraphics[width=\columnwidth]{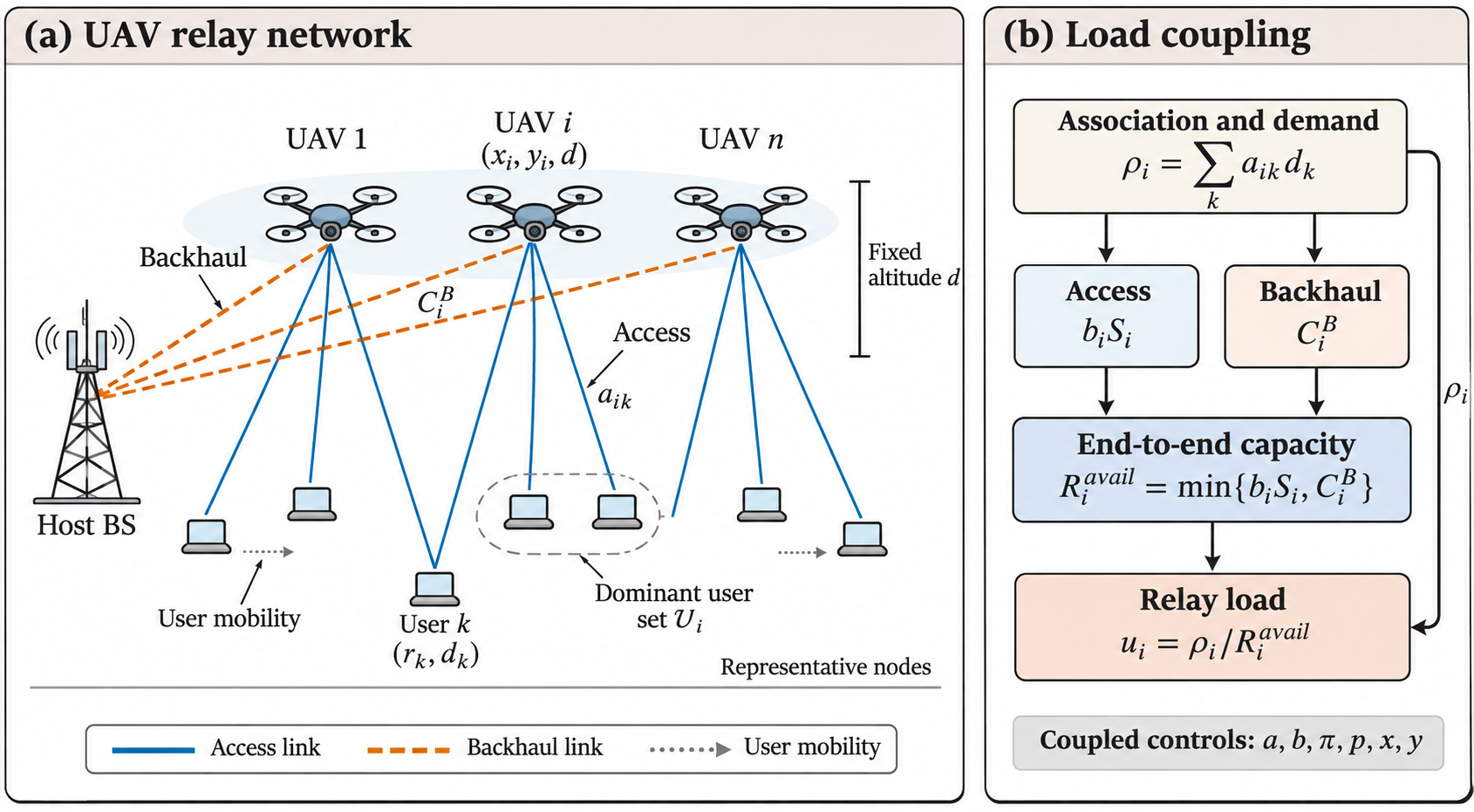}
\caption{System and load coupling. UAVs reach the host base station over the backhaul hop and serve ground users over the access hop. The load of a UAV is set jointly by association on the demand side and by position, bandwidth and power on the capacity side, and UAVs are coupled through a shared user pool and through interference.}
\label{fig:network}
\end{figure}

Consider a relay network with $n$ UAVs and $K$ users, illustrated in Fig.~\ref{fig:network}. The association weight between user $k$ and UAV $i$ is $x_{ik}\ge0$ with $\sum_i x_{ik}=1$, which gives soft association. The demand of user $k$ is $d_k>0$ and the demand carried by UAV $i$ is $\rho_i=\sum_k x_{ik}d_k$. Onboard bandwidth is split in proportion to the associated demand, $b_{ik}=x_{ik}d_k b_i/\rho_i$. The spectral efficiency $r_{ik}=\log_2(1+\mathrm{SINR}_{ik})$ is set by the channel and the power, and this paper treats $r_{ik}(X)$ as a given function of the remaining state $X$, which collects positions, power and the channel realisation, without opening its internal structure. The access rate, the end-to-end rate and the load are
\begin{equation}
\begin{split}
R_i^{\mathrm A}&=\sum_k x_{ik}b_{ik}r_{ik},\\
R_i^{\mathrm{avail}}&=\min\{R_i^{\mathrm A},\,C_i^{\mathrm B}\},\\
u_i&=\frac{\rho_i}{R_i^{\mathrm{avail}}},
\end{split}
\label{eq:model}
\end{equation}
where $C_i^{\mathrm B}>0$ is a deterministic function of position, the backhaul capacity, computed from the backhaul link budget and the bandwidth share of the UAV. Bandwidth obeys $\sum_i b_i\le B$ with $b_i>0$, the price lies in $\pi_i\in[-3,3]$, the power multiplier lies in $\tilde p_i\in[0.5,2]$, and the position lies in the field box $\mathcal F$. Notation is collected in Table~\ref{tab:notation}. The model assumptions that the analysis carries are gathered in Appendix~\ref{app:assumptions}: the demand-proportional split inside a UAV, the access-limited region in which the closed form is derived, the logarithmic knee view, the damping of the association block, and the two bounds that the implementation enforces on gain and action. Each proposition names the ones it needs, and Appendix~\ref{app:orbits} carries the arithmetic of the one-dimensional statement by hand at unit single-step bound.

\begin{table}[!t]
\centering
\caption{Notation.}
\label{tab:notation}
\footnotesize
\setlength{\tabcolsep}{4pt}
\begin{tabular}{@{}>{\raggedright\arraybackslash}p{0.32\columnwidth}>{\raggedright\arraybackslash}p{0.64\columnwidth}@{}}
\toprule
Symbol & Meaning \\
\midrule
$n,\ K$ & number of UAVs and of users \\
$x_{ik}$ & association weight of user $k$ at UAV $i$, $\sum_i x_{ik}=1$ \\
$d_k,\ \rho_i$ & user demand and demand carried by UAV $i$ \\
$b_i,\ b_{ik}$ & bandwidth of UAV $i$ and share given to user $k$ \\
$r_{ik},\ R_i^{\mathrm A}$ & spectral efficiency and access rate \\
$C_i^{\mathrm B},\ R_i^{\mathrm{avail}}$ & backhaul capacity and end-to-end rate \\
$u_i,\ \hat u_i$ & load and load view, $\hat u=m(u)$ \\
$S_i,\ c_i,\ c_i^{\mathrm B}$ & spectral sum, demand scale, backhaul load floor \\
$\pi_i,\ \tilde p_i$ & dual price and power multiplier \\
$B,\ \mathcal F$ & total bandwidth and field box \\
$\tau_u,\ \alpha,\ \kappa,\ \varepsilon$ & association damping and weights, 5, 1.0, 0.5, 0.1 \\
$u^\ast$ & load anchor, 0.8 \\
$m(\cdot),\ \Phi$ & logarithmic knee view and complete-graph equal-weight potential \\
$\eta,\ \gamma,\ h$ & executor gain, equivalent gain, physical single-step bound \\
$\beta$ & ill-conditioning gain, $\max_i\{1/S_i^2,\ \kappa/(\rho_i+\varepsilon)^2\}$ \\
\bottomrule
\end{tabular}
\end{table}

\subsection{Load View and Potential}

The load view is the logarithmic knee map
\begin{equation}
\begin{split}
m(u)&=\begin{cases}u, & u\le 3,\\[2pt] 3+3\ln\!\bigl(1+\tfrac{u-3}{3}\bigr), & u>3,\end{cases}\\
m'(u)&=\begin{cases}1,& u\le3,\\[2pt] \dfrac{1}{1+(u-3)/3},& u>3.\end{cases}
\end{split}
\label{eq:view}
\end{equation}
The map is continuously differentiable with $0<m'\le1$, so $m$ is a strictly increasing $1$-Lipschitz bijection from $[0,\infty)$ onto itself, and $m'(u)\to0$ as $u\to\infty$, which is the logarithmic compression of the overload region. The potential takes the equal-weight complete-graph form
\begin{equation}
\Phi(\hat u)=\frac12\sum_{i<j}(\hat u_i-\hat u_j)^2=\frac n2\sum_i(\hat u_i-\bar{\hat u}\,)^2,\qquad \hat u_i=m(u_i).
\label{eq:phi}
\end{equation}

\begin{proposition}[Order preservation and zero set]\label{prop:bijection}
The map $m$ is strictly increasing and continuous, so $\Phi(\hat u)=0$ holds if and only if $\hat u_1=\cdots=\hat u_n$, which holds if and only if $u_1=\cdots=u_n$. The equal-load set and the zero set of the potential therefore coincide. A clipped map $\mathrm{clip}(\cdot,0,c)$ lacks the property, since every preimage in the clipped region collapses to one point and the gradient of the potential over deep-overload differences vanishes identically.
\end{proposition}

\begin{pf}
A strictly increasing bijection preserves equality and its inverse gives both directions. The clipped map equals $c$ on $(c,\infty)$, so for $u_i\ne u_j>c$ one gets $\hat u_i=\hat u_j$ while $u_i\ne u_j$, and the gradient of $\Phi$ with respect to any change of the pair $(u_i,u_j)$ is then zero.
\end{pf}

\noindent The view also delivers the gradient chain of the potential. Using $\hat u_i=m(u_i)$ and $u_i=c_i/b_i$ in the access-limited region of Lemma~\ref{lem:homogeneous},
\begin{equation}
\frac{\partial\Phi}{\partial b_i}=\Bigl(\sum_{j\ne i}(\hat u_i-\hat u_j)\Bigr)\,m'(u_i)\,\Bigl(-\frac{u_i}{b_i}\Bigr).
\label{eq:gradchain}
\end{equation}
The gradient is strongest in the healthy band, where $m'=1$, and is compressed logarithmically but never extinguished in the deep-overload region, where $m'\to0$. That behaviour is the advantage of the view over clipping, since the recovery direction always exists and only the step conversion shrinks with depth. Table~\ref{tab:views} contrasts the three candidate views. Both clipped views have a vanishing derivative in the clipped region, which loses the order of deep overload and kills the gradient, whereas the logarithmic knee keeps a strictly positive derivative everywhere.

\begin{table}[!t]
\centering
\caption{Three load views.}
\label{tab:views}
\scriptsize
\setlength{\tabcolsep}{3pt}
\begin{tabular}{@{}lllll@{}}
\toprule
View & Definition & Derivative & Order & Failure \\
\midrule
A, clip 3 & $\mathrm{clip}(u,0,3)$ & $0$ & no & all $u>3$ \\
B, clip 30 & $\mathrm{clip}(u,0,30)$ & $0$ & no & all $u>30$ \\
C, knee & (\ref{eq:view}) & $1/\bigl(1+\frac{u-3}{3}\bigr)>0$ & yes & none, no ceiling \\
\bottomrule
\end{tabular}
\end{table}

\subsection{Why the Ratio Is the State Variable}

The state of a UAV could be described in many ways: the vector of its assigned users, the vector of its per-user rates, the bandwidth it holds, or the demand it carries. The analysis uses none of these directly and works with the ratio of carried demand to available rate, which deserves a justification and not a mere convention.

Three arguments support the choice. The first is that the ratio is the quantity that the objective penalises. The reported indicator, the Jain index, is a function of the ratio and of nothing else, and the potential is a function of the ratio as well, so a state description that already contains the argument of both is the minimal one. The second is that the ratio absorbs a whole class of changes that leave the behaviour invariant. Doubling every demand and doubling every bandwidth leaves every ratio unchanged, and the fleet behaves identically; the same holds for doubling the field and doubling the distances, up to the constant terms of the path loss. A state variable that is invariant under the transformations that the dynamics respects prevents the analysis from having to re-derive the same statement under each rescaling.

The third argument is the decisive one for the closed form. In the access-limited region the load is exactly the demand scale divided by the bandwidth, so the ratio is linear in the inverse of the decision. A statement about the load is therefore a statement about the bandwidth up to a reciprocal, and a target load is a target bandwidth. No other state description has that property: a statement about the per-user rates would carry the assignment and the spectral efficiencies, and a statement about the bandwidth alone would carry the demand scales implicitly.

The price of the choice is that the ratio is undefined when the available rate vanishes, which the implementation handles with a floor and which the analysis treats as the boundary of the model, not as a case to be handled. The floor is a numerical device and the statements of Section~III are made for positive rates. One consequence is worth naming: because the ratio is undefined at zero rate, no continuity argument can be made across the degenerate point, which is why the audit reports the working states and not the saturated ones.

\subsection{Association Block}

\begin{figure}[!t]
\centering
\includegraphics[width=\columnwidth]{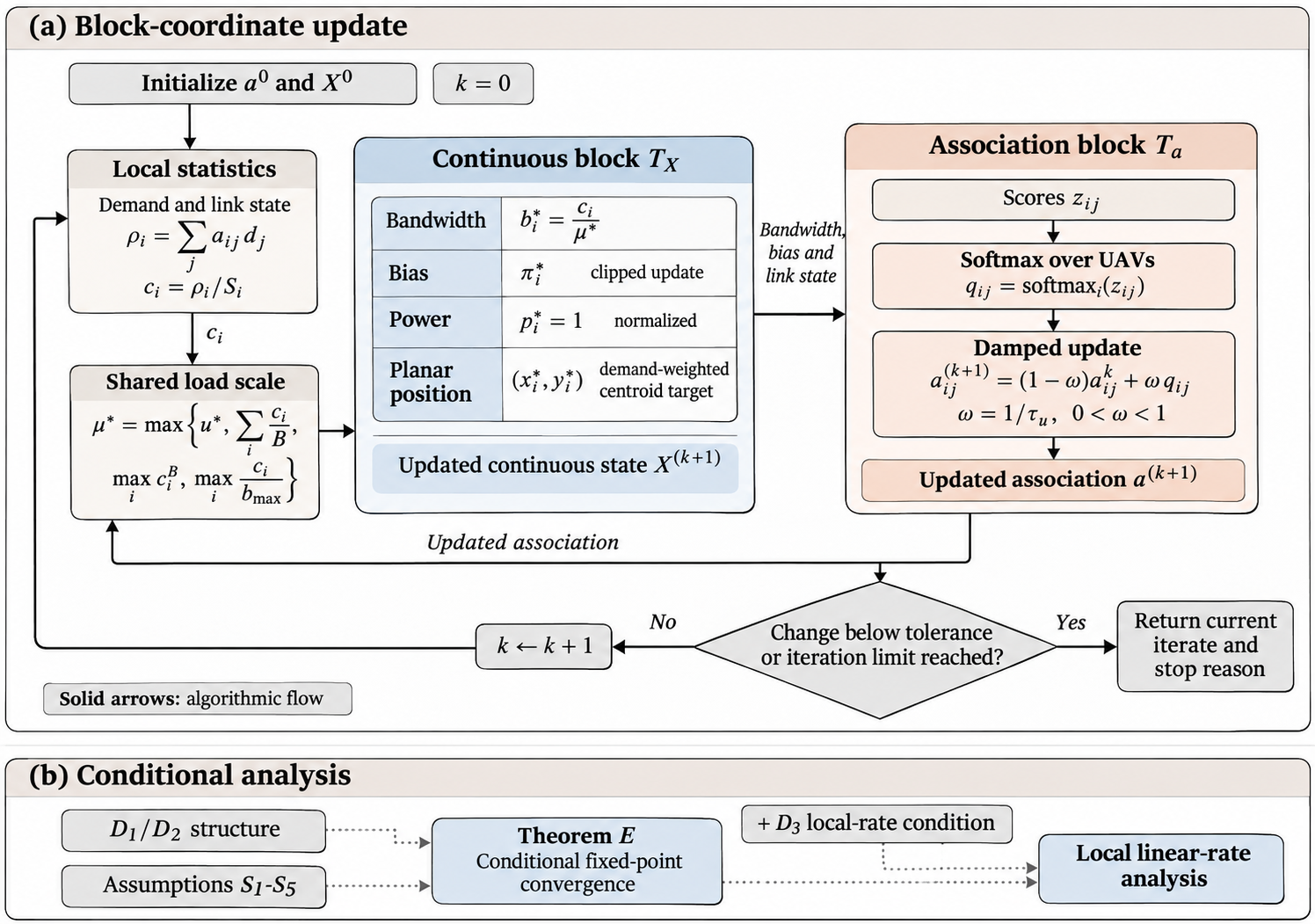}
\caption{Alternating update of association and continuous quantities. Every step first updates the association by the damped softmax of (\ref{eq:assocT}), a convex combination of the old weights and the new best response, and then recomputes the closed-form targets of bandwidth, price, power and position under the updated association. The executor composes the action from the gains.}
\label{fig:block}
\end{figure}

The association evolves by a damped softmax, drawn in Fig.~\ref{fig:block},
\begin{equation}
\begin{split}
T_x(x,X)&=\Bigl(1-\frac1{\tau_u}\Bigr)x+\frac1{\tau_u}\operatorname{softmax}\bigl(z(X)\bigr),\\
z_{ik}&=\alpha r_{ik}+\kappa\frac{R_i^{\mathrm{avail}}}{\rho_i+\varepsilon}+\pi_i,
\end{split}
\label{eq:assocT}
\end{equation}
where $X$ denotes the remaining state, so that $z$ is given, and $\tau_u=5$.

\begin{proposition}[Existence of an association fixed point]\label{prop:brouwer}
For every fixed remaining state $X$ the map $T_x(\cdot,X)$ is continuous from the product of simplices into itself, and the Brouwer fixed point theorem gives a fixed point of the association block.
\end{proposition}

\begin{pf}
The softmax is continuous and its output lies in the closure of the open simplex, so $(1-1/\tau_u)x+(1/\tau_u)\operatorname{softmax}(z)$ is a continuous convex-combination map of the closed simplex into itself. A finite product of simplices is a nonempty closed convex compact set, which is exactly the setting of the Brouwer theorem.
\end{pf}

\noindent\textbf{Remark on design meaning.} The audit shows association weights saturating deeply in steady state, since the median largest weight over the sixty measured states is of order $1-10^{-8}$, and the update map of the saturated region amplifies perturbations in the Euclidean metric, as the spectral observation of Section~V reports. The runtime guarantee of the system holds under conditions strictly weaker than uniqueness of the fixed point. The damping inertia with $\tau_u=5$ locks the association into a slow variable, the simplex bounds every weight, the damping keeps those weights slow, and the outer mechanism of Section~IV supplies the guarantee for any association configuration. The analysis therefore covers all configurations of the block, which is the strength that the layering delivers.

\subsection{Position Potential}

The position target of the direction decider is the demand-weighted centroid of the served users, and the position potential supplies the optimality of that choice.

\begin{proposition}[Unique minimiser of the position potential]\label{prop:centroid}
Define the position potential of UAV $i$ as $\Psi_i(x)=\sum_{k\in U_i}d_k\|x-r_k\|^2$ over $x\in\mathbb R^2$, with $U_i$ the set of users it serves. The potential is strictly convex and its unique global minimiser is the demand-weighted centroid
\begin{equation}
x_i^\ast=\frac{\sum_{k\in U_i}d_k r_k}{\sum_{k\in U_i}d_k}.
\end{equation}
After projection onto the field box $\mathcal F$ the point $\mathrm{clip}_{\mathcal F}(x_i^\ast)$ is the unique minimiser of $\Psi_i$ over $\mathcal F$, since the isotropic quadratic separates into independent coordinates and coordinate-wise clipping minimises each term.
\end{proposition}

\begin{pf}
The Hessian $\nabla^2\Psi_i=2\bigl(\sum_{k}d_k\bigr)I$ is positive definite, so the potential is strictly convex, and the first-order condition $2\sum_k d_k(x-r_k)=0$ has the unique solution $x_i^\ast$. On a rectangular field box the objective is a positive sum of squared coordinate differences. Clipping each coordinate of $x_i^\ast$ to its feasible interval minimises the corresponding one-dimensional quadratic, and the sum of these minima gives the unique constrained minimiser.
\end{pf}

\noindent Combined with Lemma~\ref{lem:homogeneous}, both the bandwidth direction and the position direction now carry proofs. Homogeneity gives the bandwidth target in Section~III and centroid optimality gives the position target. The price dimension has no within-snapshot optimal solution, since the dual equilibrium requires iteration, and its direction is a subgradient sign. The direction is thus three closed-form dimensions plus one subgradient dimension.

\section{Homogeneous Structure and Anchored Closed Form}

\begin{lemma}[First-order homogeneous decomposition]\label{lem:homogeneous}
Under the demand-proportional bandwidth allocation of (\ref{eq:model}),
\begin{equation}
R_i^{\mathrm A}=b_i S_i,\qquad
S_i=\sum_k\frac{x_{ik}^2 d_k}{\rho_i}\,r_{ik},
\label{eq:Si}
\end{equation}
where $S_i$ depends on association, positions, power and the channel realisation through $r_{ik}$, but not on $b_i$. In the access-limited region, where $R_i^{\mathrm A}\le C_i^{\mathrm B}$, the load is therefore
\begin{equation}
u_i=\frac{c_i}{b_i},\qquad c_i=\frac{\rho_i}{S_i},
\label{eq:ci}
\end{equation}
with a demand scale $c_i$ independent of $b_i$. In the backhaul-limited region the load is $u_i=c_i^{\mathrm B}=\rho_i/C_i^{\mathrm B}$, again independent of $b_i$, and the two branches join continuously at $b_iS_i=C_i^{\mathrm B}$, equivalently $b_i=C_i^{\mathrm B}/S_i$. The demand scale equals $c_i=u_ib_i$ numerically in the access-limited region and is read directly from the current state of the UAV.
\end{lemma}

\begin{pf}
Substituting $b_{ik}=x_{ik}d_k b_i/\rho_i$ gives $R_i^{\mathrm A}=\sum_k x_{ik}\cdot\frac{x_{ik}d_k b_i}{\rho_i}\cdot r_{ik}=b_i\sum_k\frac{x_{ik}^2d_k}{\rho_i}r_{ik}$, whose sum carries no $b_i$. The two branches of the load are $\rho_i/(b_iS_i)=c_i/b_i$ and $\rho_i/C_i^{\mathrm B}$. At the junction $b_iS_i=C_i^{\mathrm B}$ the two expressions agree, and the load as a function of $b_i$ is continuous, strictly decreasing on each branch and not differentiable at the junction, where the left derivative $-c_i/b_i^2$ meets the right derivative $0$.
\end{pf}

\begin{proposition}[Anchored closed form and feasibility bound]\label{prop:anchor}
Fix the association, positions, power and channel realisation, which fixes the vectors $c=(c_i)$ and $c^{\mathrm B}=(c_i^{\mathrm B})$. The equal-load point $u_i\equiv\mu$ is reachable through bandwidth if and only if
\begin{equation}
\mu\ \ge\ \mu_{\mathrm{low}}=\max\Bigl\{\frac{\sum_i c_i}{B},\ \max_i c_i^{\mathrm B}\Bigr\},
\label{eq:mulow}
\end{equation}
and a minimal-bandwidth allocation that attains it is $b_i=c_i/\mu$. This allocation is unique whenever $c_i^{\mathrm B}<\mu$ for every $i$. The system adopts $\mu^\ast=\max\{u^\ast,\ \mu_{\mathrm{low}}\}$ with anchor $u^\ast=0.8$. When the budget is sufficient, every UAV anchors at $u^\ast$. When the budget determines the common level and no higher backhaul floor binds, $\mu^\ast=\sum_ic_i/B$ and the closed form degenerates into the proportional split $b_i=Bc_i/\sum_jc_j$.
\end{proposition}

\begin{pf}
Equal load requires $\max\{c_i/b_i,\,c_i^{\mathrm B}\}=\mu$ for every $i$, which implies $b_i\ge c_i/\mu$ and $\mu\ge c_i^{\mathrm B}$. Minimising total bandwidth at fixed $\mu$ therefore gives $b_i=c_i/\mu$, and the budget constraint requires $\sum_ic_i/\mu\le B$. Conversely, this allocation attains load $\mu$ on every coordinate whenever $\mu$ satisfies both lower bounds in (\ref{eq:mulow}), which proves reachability and minimality. If $c_i^{\mathrm B}<\mu$, equality of the two-branch load forces $b_i=c_i/\mu$. If $c_i^{\mathrm B}=\mu$, any larger feasible bandwidth also leaves the coordinate on its backhaul branch, so uniqueness need not hold at that boundary. The anchored and budget-determined branches follow by substituting $\mu^\ast$ into the minimal-bandwidth allocation.
\end{pf}

\noindent\textbf{Corollary on the bandwidth direction and on communication.} The bandwidth target $b_i^\ast=c_i/\mu^\ast$ of the direction decider depends only on the local quantity $c_i$ and on the broadcast sum of demand scales $\sum_jc_j$, together with the public budget $B$; the backhaul caps enter the feasibility bound $\mu_{\mathrm{low}}$ of Proposition~\ref{prop:anchor} and not the computation of $b_i^\ast$. The broadcast mean price enters the power update and not the bandwidth target. The demand scale $c_i=u_ib_i$ is read directly from the current state, which Lemma~\ref{lem:homogeneous} guarantees to be independent of $b_i$, and the audited drift is $2.7\times10^{-16}$. No anchoring strength appears in the implementation, because hard anchoring takes the zero directly and an anchoring weight would change only the curvature away from the zero and not the zero itself.

\noindent\textbf{Discussion.} The two regions of Proposition~\ref{prop:anchor} each have an economic reading. In the anchored region bandwidth is the abundant resource, every UAV takes what its demand scale asks for, the load is pressed uniformly onto the anchor and fairness is maximal. In the budget region total bandwidth cannot press every UAV onto the anchor, the optimal response is to spend the budget in proportion to demand scale, and the common load $\mu=\sum c/B$ then sits above the anchor, which places the system at the maximal equilibrium permitted by the budget. The backhaul floor $c_i^{\mathrm B}$ is a per-UAV ceiling. When the backhaul capacity of a UAV is low enough that $\rho_i/C_i^{\mathrm B}>\mu$, that UAV is pinned on the backhaul branch and extra bandwidth does nothing. The situation does not arise under the parameter setting of \cite{r69}, where the backhaul is normally ample, yet it is kept in the structure so that the model stays semantically complete for backhaul-tight scenarios such as distant UAVs or large backhaul delay.

\subsection{A Small Instance Worked Through}

The closed form is short enough to be checked by hand, and one instance fixes the reading of the two regions. Take four UAVs with demand scales $c=(2,3,5,10)$ in the same unit as the budget and a budget $B=25$. The budget-determined level is $\sum_jc_j/B=0.8$, which equals the anchor $u^\ast=0.8$, so the two regions meet at this budget: the target is $b^\ast=(2.5,3.75,6.25,12.5)$, every UAV sits at $u=0.8$, the total bandwidth spent is $25$ and the Jain index is one.

Lower the budget to $B=20$ and the anchor is no longer reachable, since $\mu_{\mathrm{low}}=\sum_jc_j/B=1>u^\ast$. The level is then $\mu^\ast=1$ and the target is $b^\ast=(2,3,5,10)$, which again spends the whole budget, again leaves every UAV at a common load and again gives a Jain index of one. The two cases differ in the level of the common load and not in its uniformity, so the anchor governs the level only while the budget is large enough to afford it.

Now move the crowd so that the third demand scale doubles to $10$. The sum becomes $25$ and the budget branch gives $u^\ast=25/20=1.25$ together with $b^\ast=(1.6,2.4,8,8)$, which again sums to the budget and again gives identical loads. Two things are worth reading off the numbers. The UAV whose demand did not change still loses bandwidth, from $2$ to $1.6$, because the budget is finite and the extra demand of another UAV is paid by every UAV in proportion to its own scale. And the common load rises from one to $1.25$, which is exactly the ratio of the new sum to the budget, so the price level and the load are the same object in the budget-bound region.

The access-limited branch appears when the budget is larger than the anchor requires. With the original demands and $B=40$ the level stays at the anchor, so the target is still $b^\ast=(2.5,3.75,6.25,12.5)$, which sums to $25$ and leaves $15$ of the budget unused. That is why the total bandwidth constraint is rarely binding in a lightly loaded fleet, and the branch matters in the opposite case, a fleet with a small budget and a large backhaul ceiling, where the access hop pins the load while the budget stays slack. The direction decider does not need to know which branch is active: it computes $b_i^\ast=c_i/\mu^\ast$ with $\mu^\ast=\max\{\sum_jc_j/B,\ \max_ic_i^{\mathrm B}\}$ and then rescales to the budget, which reproduces both branches and the boundary between them.

\subsection{Potential and Fairness Index}

The balance indicator that the system reports is the Jain fairness index defined on the raw load \cite{r37}, $J(u)=\bar u^2/(\bar u^2+\operatorname{var}_u)$. Its relation to the potential is as follows.

\begin{proposition}[Potential versus fairness index]\label{prop:jain}
Write $\operatorname{var}_u=(1/n)\sum_i(u_i-\bar u)^2$ and $\operatorname{var}_{\hat u}=(1/n)\sum_i(\hat u_i-\bar{\hat u})^2$. First, the zero sets coincide, since $\Phi=0$ holds if and only if $\operatorname{var}_{\hat u}=0$, if and only if $\operatorname{var}_u=0$, if and only if $J(u)=1$. Second, on the linear band of the view, where $m(u)=u$, the two objects are monotone functions of the same variance at a fixed mean, since $\Phi=(n^2/2)\operatorname{var}_u$ and $J(u)=1/(1+\operatorname{var}_u/\bar u^2)$, so at fixed $\bar u$ any change that lowers $\Phi$ raises $J$. Third, outside that band the view preserves the order of the loads and the sign of the dispersion but not its magnitude, so the two objects share the optimal direction and the zero set without being interchangeable as measures of distance.
\end{proposition}

\begin{pf}
For the first claim, Proposition~\ref{prop:bijection} makes $\operatorname{var}_{\hat u}=0$ equivalent to equality of all $\hat u$, hence to equality of all $u$ and to $\operatorname{var}_u=0$, and $J(u)=1$ is equivalent to $\operatorname{var}_u=0$ because $\bar u>0$. For the second, the derivative of $J=1/(1+\operatorname{var}_u/\bar u^2)$ with respect to $\operatorname{var}_u$ is negative at fixed positive $\bar u$, and $\Phi=\frac n2\sum(\hat u_i-\bar{\hat u})^2=\frac{n^2}{2}\operatorname{var}_{\hat u}$ is strictly increasing in $\operatorname{var}_{\hat u}$, which gives the monotone pair on the linear band. For the third, $m$ is strictly increasing by Proposition~\ref{prop:bijection}, so equal viewed loads are equal raw loads, while the curvature of $m$ rescales the dispersion by a load-dependent factor, which is why the magnitudes of $\Phi$ and $\operatorname{var}_u$ cannot be exchanged outside the linear band.
\end{pf}

\noindent The proposition separates the roles of the two quantities. The potential, which is the variance of the viewed load, is the object being optimised, and the view places deep overload on one scale while preserving order. The fairness index, the Jain form on the raw load, is the object being reported, with a scale that is intuitive and directly comparable with the literature. The zeros coincide and the directions agree, so that optimising the potential cannot worsen fairness inside the linear band, and outside that band the two objects still share the optimum and the sign of the dispersion.

\subsection{Why a Global Contraction Cannot Serve as the Certificate}

A common way to certify a decentralised resource allocation is to write the update as a contraction in a global metric and to quote the contraction rate. That route is closed here, and it is worth saying exactly why, because the reason shapes what this paper proves.

The map that acts on the fleet is a composition of five nonlinearities: the softmax over users inside the association block, the clamp that turns an access rate into a load, the maximum that selects the serving hop, the proportional rescaling that keeps the shares on the budget simplex, and the trunk of the executor network with its bounded gain. Each of the five is individually well behaved. The softmax is smooth and its Jacobian has spectral radius bounded by the largest class probability. The clamp has a Jacobian of rank one, the maximum is piecewise linear with unit slope, the rescaling is a projection onto a simplex, and the network is a bounded smooth function with a derivative bounded by the product of its weight norms. A product of five such maps has no reason to contract: the softmax Jacobian can be small while the rate-to-load map expands, the projection can expand along the tangent of the simplex, and the network derivative can compensate for the others at any point of the state space.

The numerical evidence in Section~V-B makes the failure concrete. The one-step map of the association block, measured on states of the audit, has a spectral radius between $1.0\times10^{2}$ and $3.1\times10^{3}$ depending on the operating point, so no single contraction factor exists that covers the fleet across its regimes. The reason is structural: the load is a ratio of a carried demand to an available rate, and the available rate is itself a function of how the users are assigned, so the map contains an inverse of its own output in the saturated regime. The softmax curvature does not save the certificate, because at saturation the class probabilities are nearly uniform in the logit and the curvature is tiny while the ratio still moves.

What replaces the certificate is a decomposition along the layers, and each layer has a property that holds for its own reasons. The direction is homogeneous and anchored, so its fixed point is available as a closed form, with no contraction needed to reach it. The executor is bounded, so its output never leaves the physical envelope regardless of the state. The one-dimensional gain iteration has a complete classification, with geometric convergence below the critical gain and a bounded two-state family above it. The tracking bound follows from an explicit recursion with quantified disturbances. The spectral statement is then a statement about how far one block may amplify, with the amplification measured, not assumed, and the conclusion is that the composite is bounded and that its oscillations are confined to the conservative gain range of a single dimension.

The certificate is not a global contraction, and it does not have to be one. A claim of contraction would have to be either restricted to a neighbourhood that the experiments leave, or paid for with a conservative rate that predicts divergence where the fleet in fact converges. The layered statement is the stronger object of the two where it matters: it stays valid across the whole gain range and states the price of the gain explicitly, below two convergence, above two a bounded family of cycles whose amplitude is known.

\section{Executor Dynamics and Step-Size Bounds}

The executor runs a truncated gain iteration independently on every dimension, and the iteration admits a complete closed form in one dimension.

\subsection{One-Dimensional Truncated Gain Dynamics}

Fix one dimension with target $x^\ast$, current value $x_t$, error $e_t=x_t-x^\ast$, physical single-step bound $h>0$ and equivalent gain $\gamma>0$. The iteration is
\begin{equation}
e_{t+1}=e_t-\mathrm{clip}(\gamma e_t,\ \pm h).
\label{eq:1d}
\end{equation}

\begin{theorem}[Truncated gain iteration]\label{thm:1d}
For (\ref{eq:1d}) the following hold. First, when $0<\gamma<2$ the error vanishes: while $|e_t|>h$ a saturated step removes exactly $h$ from the magnitude, a saturated step from $h/\gamma<|e_t|\le h$ sends the magnitude to $h-|e_t|\le h-h/\gamma$, and once $|e_t|\le h/\gamma$ the linear step multiplies the magnitude by $|1-\gamma|<1$, so the error decays geometrically and at $\gamma=1$ the linear step reaches the target in one step. Second, at $\gamma=2$ every orbit either reaches the target in finitely many steps or enters the pair $\{a,-a\}$, where $0<a\le h/2$ is fixed by the initial value, so the error does not tend to zero. Third, when $\gamma>2$ each $p\in[h/\gamma,\,(1-1/\gamma)h]$ generates a two-state cycle $\{p,\,p-h\}$, these are the only period-two orbits, and their amplitudes are at most $(1-1/\gamma)h<h$. An orbit reaches the target exactly when some iterate equals $\pm h$, which happens at $e_0=h$ and at every odd multiple of $h$; every other orbit ends in one of these cycles, so the family is the generic long-run behaviour and the error of a generic initial value does not tend to zero. Fourth, for every $\gamma>0$ and every initial value, $|e_t|\le\max\{|e_0|,\,h\}$ holds at all $t$, the slab $|e|\le h$ is positively invariant, and the error never diverges.
\end{theorem}

\begin{pf}
For the first claim, consider the magnitude $V=|e_t|$. If $V>h$ the clip saturates and $e_{t+1}=e_t-h\,\mathrm{sgn}\,e_t$, so $V$ decreases by $h$ exactly. If $h/\gamma<V\le h$ the clip also saturates and $V$ maps to $h-V$; since $V>h/\gamma$ gives $h-V<h-h/\gamma$, and $h-h/\gamma\le h/\gamma$ because $\gamma\le2$, that image lies in the linear range, so the magnitude strictly decreases and the next step is linear. Once $V\le h/\gamma$ the linear step gives $|e_{t+1}|=|1-\gamma||e_t|$ with $|1-\gamma|<1$ and the error decays geometrically to zero; at $\gamma=1$ the linear step maps $e$ to $0$ directly.

For the second claim, at $\gamma=2$ the same two saturated cases send $V>h$ to $V-h$ and $h/2<V\le h$ to $h-V$, both of which are at most $h/2$ or exactly zero, so after finitely many steps either the orbit hits the target, which is fixed, or $V\le h/2$. From any $V\le h/2$ the step is linear with $e_{t+1}=-e_t$, so the orbit alternates between $+a$ and $-a$ with $a=V$ fixed by the initial value, and $a\le h/2$.

For the third claim, let $p\in[h/\gamma,(1-1/\gamma)h]$. From $p$ the gain saturates, since $\gamma p\ge h$, and the image is $p-h$; from $p-h$ one has $\gamma(h-p)\ge h$ by the upper end of the interval, so the image is $p$ again and the pair is a cycle, with amplitude $\max\{p,h-p\}\le(1-1/\gamma)h<h$. If $e_t\ne0$ and $|e_t|\le h/\gamma$ the linear step amplifies the magnitude by the factor $\gamma-1>1$ until it leaves that range and saturates again, so a bounded orbit cannot converge to the target, and the period-two cycles just described are the only ones, since a two-cycle $\{p,q\}$ with $q\ne p-h$ would need a linear image to coincide with a saturated image. Inside the slab $|e|\le h$ both branches map into $[-h,h]$, so every orbit stays in the slab. A saturated step moves the state by exactly $h$ and a linear step sends $e$ to $(1-\gamma)e$, so the target is met only when an iterate equals $\pm h$ or $0$; the odd multiples of $h$ are such initial values, and every orbit that misses this set circulates in the slab without converging.

For the fourth claim, argue by induction on the magnitude. A saturated step from $V>h$ gives $V-h\le V$. A saturated step from $h/\gamma<V\le h$ gives $h-V<h\le\max\{V,h\}$. A linear step gives $(\gamma-1)V\le(\gamma-1)h/\gamma<h$. Hence the magnitude never exceeds $\max\{|e_0|,h\}$, and once it lies in the slab it stays there.
\end{pf}

\begin{corollary}[Step count to a target accuracy]\label{cor:steps}
For $0<\gamma<2$, $\gamma\ne1$ and $0<\varepsilon<h/\gamma$ the number of steps to reach accuracy $\varepsilon$ is at most
\begin{equation}
T\ \le\ \Bigl\lceil\frac{\max\{|e_0|-h/\gamma,\,0\}}{h}\Bigr\rceil+1+\Bigl\lceil\frac{\ln(\varepsilon/(h/\gamma))}{\ln|1-\gamma|}\Bigr\rceil,
\end{equation}
where the first term counts the saturated phase, which removes $h$ per step, and the last two terms count the linear phase, which starts once $|e|\le h/\gamma$ and then decays geometrically. At $\gamma=1$ the target is reached exactly within $\lceil\max\{|e_0|-h,0\}/h\rceil+1$ steps. At $\gamma=2$ the error enters the positively invariant interval $|e|\le h/2$ in finitely many steps and then either sits at the target or oscillates there, so no step count to an accuracy below the amplitude of that oscillation exists.
\end{corollary}

\noindent Theorem~\ref{thm:1d} matches the hard bound of the system executor, $\eta\in[0.50,2.01]$, which composes a multiplicative gain $\eta=\exp(0.7\tanh(\cdot))$ with a fixed relaxation of 0.2. Dimensions whose equivalent gain slightly exceeds two sit in the two-state family, with an amplitude of at most $(1-1/\gamma)h$, which is $0.503h$ at the upper end of the hard bound and therefore just over half a single step. Dimensions whose equivalent gain stays below one converge geometrically but slowly. Learning chooses the speed between those cases, and the fourth statement guarantees boundedness whatever gain is learned.

\noindent\textbf{Parallel dimensions and composite safety.} The five action dimensions of the system each run (\ref{eq:1d}) independently. The one-dimensional statements therefore parallelise directly into a per-coordinate bound, since every dimension satisfies the three statements of Theorem~\ref{thm:1d} and the distance of the state vector to the target is capped by the Euclidean combination of the per-dimension bounds. Dimensions interact through the environment, as bandwidth changes load, load changes association attractiveness and position changes the channel, yet the interaction never enters the execution law, which reads only the target and the current value. Every effect of the interaction appears as motion of the target and is absorbed by the perturbation term of Proposition~\ref{prop:track}. That division of labour, safe per dimension and coupled through the perturbation, is why five independent one-dimensional laws replace a five-dimensional controller. The fixed relaxation of 0.2 on the bandwidth and price dimensions presses the equivalent gain to $\gamma=\eta\times0.2\le0.40$, inside the convergent range of statement one, so those two dimensions keep geometric convergence even at the upper end of the hard bound. The position dimension has no relaxation, with $\gamma=\eta\le2.01$, which just crosses the threshold of two into the two-state family, and the resulting amplitude of at most $0.503h$ appears physically as a slight adjustment near the target.

\subsection{Perturbed Tracking Bound}

The environment is non-stationary, through crowd motion and shadow refresh, and information is limited, through centroid sampling lag and stale global quantities, so both the target $x^\ast$ and the update carry perturbations.

\begin{proposition}[Tracking error ball]\label{prop:track}
Suppose a one-dimensional closed-loop map is $\rho$-Lipschitz in the error with $\rho<1$ and that every step adds a perturbation of magnitude at most $\delta_{\max}$. Then
\begin{equation}
|e_t|\ \le\ \rho^t|e_0|+\frac{\delta_{\max}}{1-\rho},
\end{equation}
so the error trajectory enters any ball of radius $\delta_{\max}/(1-\rho)$ enlarged by an arbitrarily small margin, and it does not diverge in time.
\end{proposition}

\begin{pf}
The recursion $|e_{t+1}|\le\rho|e_t|+\delta_{\max}$ is unrolled step by step, which gives the bound for every $t$ and every initial value. For the stochastic form, suppose the perturbations of different steps are mutually independent with zero mean and variance at most $\sigma^2$. The error is then a sum of $t$ geometric weights applied to independent zero-mean terms, so its variance is at most $\sigma^2\sum_{k<t}\rho^{2k}=\sigma^2(1-\rho^{2t})/(1-\rho^2)$, and a sub-Gaussian bound on that sum gives, for any $\delta\in(0,1)$,
\begin{equation}
|e_t|\ \le\ \rho^t|e_0|+\sqrt{2\ln(2/\delta)}\;\sigma\sqrt{\frac{1-\rho^{2t}}{1-\rho^2}}
\end{equation}
with probability at least $1-\delta$. The deterministic form holds when the bound $\delta_{\max}$ is enforced on every step, and the high-probability form holds under the independence and zero-mean premises on the noise, which are premises about the noise and not about the map.
\end{pf}

\noindent Table~\ref{tab:disturb} lists each perturbation source and the term it contributes. Target drift is capped by the crowd speed limit, which bounds centroid motion per step by a quantity of order $v_{\max}\tau$. Measurement noise is Gaussian with zero mean and a standard deviation at the ten-metre level, so the deterministic recursion takes its empirically bounded part at that scale and the residual tail is what the high-probability form costs. Shadow switching is a jump at the hundred-step scale, whose size follows from a path-loss difference bounded by $\sigma\le13.2$ dB. The truncation residual is the cycle amplitude of Theorem~\ref{thm:1d}. Each dimension carries its perturbations into its own tracking ball.

\begin{table}[!t]
\centering
\caption{Perturbation sources and bounds, decomposed as $\delta_{\max}$ of Proposition~\ref{prop:track}.}
\label{tab:disturb}
\scriptsize
\setlength{\tabcolsep}{3pt}
\begin{tabular}{@{}l>{\raggedright\arraybackslash}p{0.28\columnwidth}>{\raggedright\arraybackslash}p{0.38\columnwidth}@{}}
\toprule
Source & Mechanism & Bound \\
\midrule
Target drift & crowd motion moves the centroid & $\le v_{\max}$ times step, 2.78 m/s by 0.1 s \\
Measurement noise & Gaussian noise on sampled user positions & ten-metre level \\
Stale information & global quantities refreshed every 0.5 s & drift times five steps \\
Shadow switching & link shadow refreshed every hundred steps & $\le13.2$ dB path-loss difference \\
Truncation residual & physical single-step clip & cycle amplitude of Theorem~\ref{thm:1d} \\
\bottomrule
\end{tabular}
\end{table}

\subsection{Ill-Conditioning Gain and Step-Size Ceiling}

\begin{proposition}[Concentration of ill conditioning and the step bound]\label{prop:beta}
Every ill-conditioned element of the load chain concentrates in the region of small spectral sum and small demand, through
\begin{equation}
\beta=\max_i\Bigl\{\frac{1}{S_i^2},\ \frac{\kappa}{(\rho_i+\varepsilon)^2}\Bigr\}.
\end{equation}
First, the closed-loop Lipschitz constant of the bandwidth dimension obeys $L_F\le\eta\,C_\ast\,\sqrt{\beta}$, where $C_\ast$ collects the composition constants of the softmax and of the view derivative. Second, $\eta<\eta^\ast=C_\ast^{-1}\beta^{-1/2}$ is a sufficient condition for contraction, so a fixed $\eta$ cannot certify contraction over the whole range once $\beta$ grows large enough. Third, the two entries of $\beta$ respond to opposite ends of the load range, the first growing as the spectral sum falls and the second growing as the carried demand falls, while the gain from a bandwidth action to the load, $u_i^2/c_i=\rho_i/(S_ib_i^2)$, grows with the carried demand and as the spectral sum falls, so a fixed $\eta$ violates $\eta^{\ast}$ at whichever end an instance sits.
\end{proposition}

\begin{pf}
Decompose the chain. First, the sensitivity of $c_i=\rho_i/S_i$ to the environment, which enters $S$ through $r$, obeys $\|\partial c_i/\partial r\|\le\rho_i/S_i^2$, since a relative error of $S$ is amplified by the factor $1/S$, and the lower scale of $S$ is the first entry of $\beta$. Second, the derivative of the rate-margin term $R_i^{\mathrm{avail}}/(\rho_i+\varepsilon)$ in the association logit with respect to $\rho_i$ equals $-\kappa R_i^{\mathrm{avail}}/(\rho_i+\varepsilon)^2$, which is the second entry of $\beta$. The view derivative $m'\le1$ of (\ref{eq:view}) and the $\tfrac14$-Lipschitz property of the softmax introduce no new ill conditioning. The gain from a bandwidth action to the load is $\|\partial u_i/\partial b_i\|=u_i/b_i=c_i/b_i^2=u_i^2/c_i=\rho_i/(S_ib_i^2)$, which grows with the carried demand and as the spectral sum falls. Multiplying the bounds of the individual links and merging the common factor into $\sqrt\beta$ gives $L_F\le\eta C_\ast\sqrt\beta$. The contraction condition $L_F<1$ is the second claim, and the third follows directly from the two entries.
\end{pf}

\noindent Proposition~\ref{prop:beta} is the mathematical ground of a design that learns speed while the direction stays closed form. The direction comes from Propositions~\ref{prop:anchor} and \ref{prop:centroid} and does not contain $\beta$, whereas the speed must adapt to $\beta$ and can only move inside the hard bound, which Theorem~\ref{thm:1d} keeps safe even beyond the bound. The proposition also matches the experimental observation of \cite{r69}. The fixed-step executor stays close to the learned executor at low and medium pressure and collapses exactly in the deep-overload setting of a hundred UAVs under burst load, where the spectral sum falls and the carried demand rises at once, the bandwidth-to-load gain grows by orders of magnitude and a fixed $\eta$ violates $\eta^\ast$.

\subsection{Two Ends of the Speed Range}

The gain range of the implementation has two ends and each of them answers a different question, which is why the experiments use two further executors, the constant-speed one and the fixed-step one.

The lower end is the constant gain that the potential-based heuristics of the literature use. Walking the direction at a fixed distance per control period is the special case $\eta\equiv\eta_0$ of the executor, and the analysis says exactly what such a law can and cannot do. It can converge whenever the constant stays below the bound of Theorem~\ref{thm:1d} along the whole trajectory, which is the case when the fleet stays inside one regime. It cannot keep that condition when the ill-conditioning gain of Proposition~\ref{prop:beta} varies along the trajectory, since a single constant has to satisfy a bound that moves. The fixed-step executor in the experiments is that law with the constant taken at the upper end of the admissible range, which is the most favourable choice available to a constant.

The upper end is the gain at which the dimension is driven to its physical bound on every step, which is the constant-speed executor. It converges fastest when the residual is large and it is the law that the hard clip makes safe, since every step moves by exactly h and the boundedness statement does not depend on the gain. It oscillates at the critical gain and above it, which is the price of taking the largest admissible speed.

The learned executor sits between the two and its argument is narrow. What it adds is a rule that chooses the speed per dimension from the observation of the residual, the load and the direction, which is the quantity that Theorem~\ref{thm:1d} shows to be the only thing that matters in one dimension. The conclusion of the gain analysis is therefore not that learning is necessary in general, and no such claim is made. It is that the object whose choice determines survival is the speed, that the safe speed is bounded per dimension by a quantity that varies with pressure, and that a state-dependent gain is the natural way to track it.

\subsection{Gain Bound and Training}

The output range of the executor network, $\eta=\exp(0.7\tanh u)\in[0.50,2.01]$, removes any safety constraint from training, since every network output is mapped into the hard bound and the boundedness statement of Theorem~\ref{thm:1d} holds identically. TD3 training selects the speed distribution of the dimensions inside that bound, where the critic estimates the per-UAV return over a reward that aggregates balance, load and throughput, and the actor moves the logits along the direction of improved return. Zero initialisation of the output layer starts training from $\eta\equiv1$, the pure potential policy, which is an executable initial strategy, and every intermediate checkpoint corresponds to a safe executor. Safety is thus built into the action space, where the physical bound delimits the safe domain, which is stronger than adding a constraint penalty to the training objective, because the structural guarantee precedes and does not depend on the training process.

\section{Numerical Verification}

The numerical audit runs over $n=3,12,30,60,100$, the two scenarios of flow and burst, and seeds $0,1,2$, which gives thirty instances and sixty states from snapshots at steps 700 and 1500. The setting is the pure field with $\eta\equiv1$, the UMi channel class, a fixed 6 dB shadow and no physical engine. Table~\ref{tab:audit} states the setting and the scale, and the results fall into four groups.

\begin{table}[!t]
\centering
\caption{Setting and scale of the numerical audit.}
\label{tab:audit}
\footnotesize
\setlength{\tabcolsep}{4pt}
\begin{tabular}{@{}>{\raggedright\arraybackslash}p{0.30\columnwidth}>{\raggedright\arraybackslash}p{0.66\columnwidth}@{}}
\toprule
Item & Content \\
\midrule
scale, scenario, seed & $n=3,12,30,60,100$, flow and burst, $s=0,1,2$ \\
states & one snapshot at step 700 and one at 1500 per instance, sixty in total \\
setting & pure field $\eta\equiv1$, UMi channel class, fixed 6 dB shadow, no physical engine \\
diagnostics & homogeneity residual D1, energy decomposition D2, spectral estimates by ARPACK \\
results & per-cell data of thirty instances, aggregate statistics of sixty states, convergence sequences of three instances \\
\bottomrule
\end{tabular}
\end{table}

\subsection{Check of Homogeneity}

\begin{figure}[!t]
\centering
\includegraphics[width=\columnwidth]{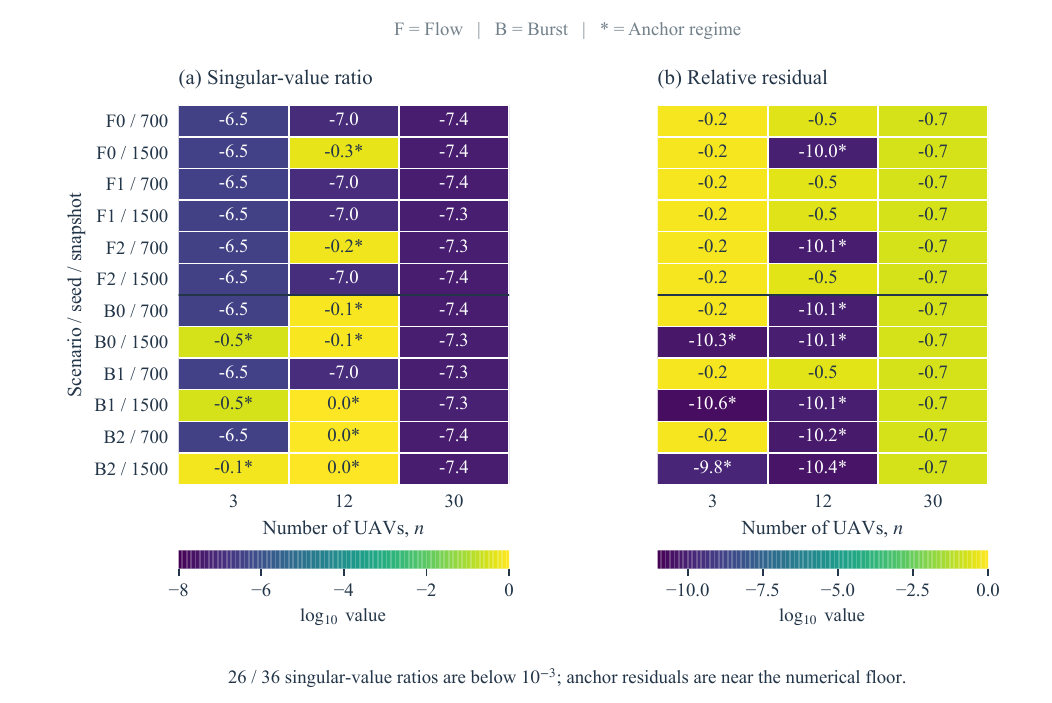}
\caption{Homogeneity residual D1 as a heat map over the operating region, with scale on the horizontal axis and state on the vertical axis and colour giving the demand-scale drift on a logarithmic scale. All sixty states lie inside the machine-precision band, and both the anchored branch and the budget branch hold.}
\label{fig:d1}
\end{figure}

The invariance of the demand scale in Lemma~\ref{lem:homogeneous} holds on all sixty states, as Fig.~\ref{fig:d1} shows. The drift of $c_i=u_ib_i$ under bandwidth perturbation is of order $10^{-16}$, that is at machine precision, and the closed-form bandwidth agrees with the feasible optimum. Homogeneity is the foundation of the entire closed-form chain, and the audit finds no deviation across two snapshots, five scales and two scenarios.

\subsection{Boundary of the Spectral Certificate}

\begin{figure}[!t]
\centering
\includegraphics[width=\columnwidth]{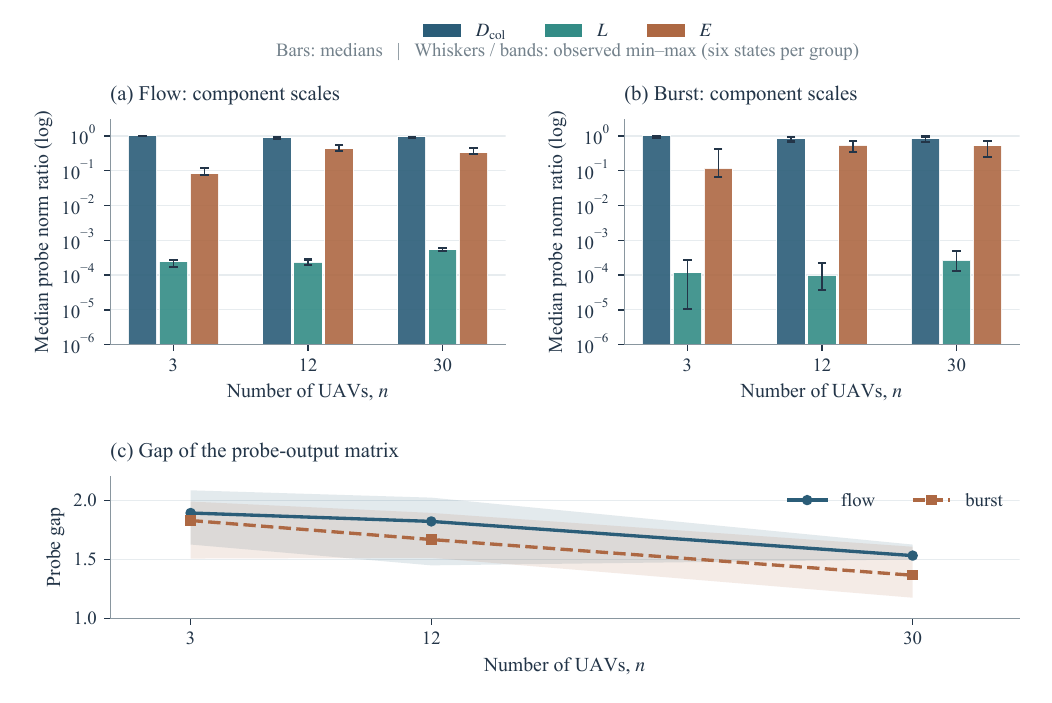}
\caption{Energy decomposition of the one-step association map, the D2 probe. Per-user independent blocks dominate, and the cross-user coupling terms are small. The decomposition keeps the same shape across the whole scale axis.}
\label{fig:d2}
\end{figure}

\begin{figure}[!t]
\centering
\includegraphics[width=\columnwidth]{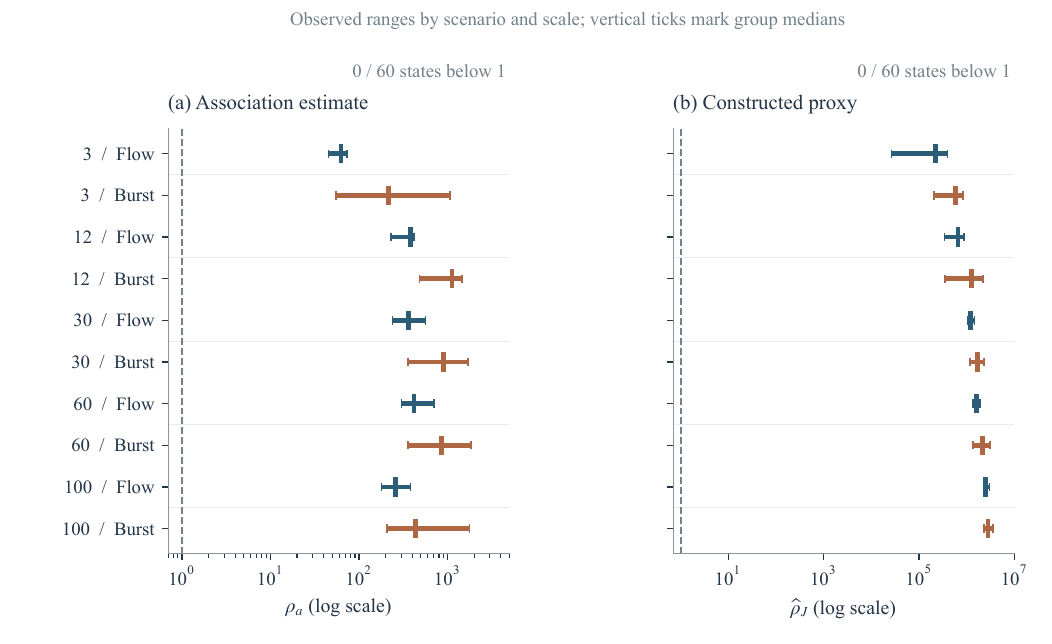}
\caption{Record of local spectral quantities, with states on the horizontal axis and both the spectral estimate and the certificate interval on a logarithmic vertical axis. The Euclidean estimate reaches $10^{1}$ to $10^{3}$ on saturated states, where the linear certificate does not hold, whereas the local estimate near the static fixed point is about one.}
\label{fig:cert}
\end{figure}

The sufficient condition for linear convergence, a closed-loop spectral radius below one, fails on all sixty measured states, as Figs.~\ref{fig:d2} and \ref{fig:cert} and Table~\ref{tab:spectral} report. The operator being measured must be stated before the numbers are read. It is the one-step association map with the whole chain recomputed from the perturbed association at every matrix-vector product, so that the bandwidth shares, the transmit power, the interference denominator and the available rate all respond to the perturbation, while the position, the price and the demand stay at their snapshot values. Its Jacobian splits into two parts, the diagonal relaxation $0.8I$ of the damped softmax and a tangential term, the product of the softmax Jacobian with the derivative of the attention logits. In the saturated regime the relaxation contributes the diagonal factor and the amplification sits in the logit derivative, since a change of association moves the interference denominator and the available-rate ratio at a rate that grows with the serving signal-to-noise ratio. The radius over the sixty states at steps 700 and 1500 runs from $50$ to $3107$, that is $10^{1.7}$ to $10^{3.5}$, on the same states whose largest association weight is of order $1-10^{-8}$. The system nevertheless converges and stays bounded on every measured state. A Euclidean derivative norm above one at saturation states the conditioning of the coordinates and not a divergence of the map, which still sends the simplex into itself, and the criterion of a radius below one presupposes an interior non-saturated point with a Euclidean metric aligned with the dynamics, while neither premise holds at a saturated boundary.

The audit records two further spectral quantities for comparison. On the static instance the ARPACK estimate in a neighbourhood of the fixed point, with difference steps from $10^{-4}$ to $10^{-10}$, is stable between $0.98$ and $1.03$, and the three snapshots at steps 3000, 6000 and 10000 give $1.027$, $0.979$ and $0.995$. The local spectral radius near the fixed point is thus about one, which corresponds to a critical rate, either slow convergence with a geometric factor approaching one or a sublinear profile. On the time-varying flow instance the trace estimate under the same setting is of order $10^{2}$, with the three snapshots at 308, 374 and 471, and that estimate takes finite differences along a moving trajectory, so it contains the drift of the target itself and is not a fixed-point spectrum. The two orders of magnitude differ precisely because spectral compression holds approximately inside the certificate boundary, near a static fixed point, and spectral quantities are unusable outside it, in the saturated and time-varying regimes, where convergence is covered by the mechanism of Section~IV. That mechanism has no certificate boundary and delivers its guarantee over the whole operating region, saturated and time-varying included.

\begin{table*}[!t]
\centering
\caption{Spectral observations, aggregated by seed median.}
\label{tab:spectral}
\footnotesize
\setlength{\tabcolsep}{6pt}
\begin{tabular}{@{}llrrrr@{}}
\toprule
$n$ & Scenario & $\bar u$ & $\rho_a$ at 700 & $\rho_a$ at 1500 & median $a_{\max}$ \\
\midrule
3 & burst & 0.845 & 56.6 & 758 & $1-4.2\mathrm{e}{-10}$ \\
3 & flow & 0.851 & 71.4 & 70.8 & $1-6.5\mathrm{e}{-10}$ \\
12 & burst & 0.866 & 386 & 1463 & $1-1.1\mathrm{e}{-9}$ \\
12 & flow & 0.873 & 308 & 375 & $1-1.9\mathrm{e}{-9}$ \\
30 & burst & 0.952 & 304 & 803 & $1-5.4\mathrm{e}{-9}$ \\
30 & flow & 0.959 & 376 & 449 & $1-7.1\mathrm{e}{-9}$ \\
60 & burst & 0.908 & 390 & 1553 & $1-2.8\mathrm{e}{-9}$ \\
60 & flow & 0.912 & 428 & 373 & $1-3.7\mathrm{e}{-9}$ \\
100 & burst & 0.982 & 271 & 1076 & $1-1.8\mathrm{e}{-8}$ \\
100 & flow & 0.996 & 172 & 270 & $1-1.8\mathrm{e}{-8}$ \\
\midrule
\multicolumn{6}{@{}l@{}}{\small static fixed point by ARPACK with $\epsilon$ swept over $10^{-4}$ to $10^{-10}$, $1.027/0.979/0.995$ at the three snapshots} \\
\multicolumn{6}{@{}l@{}}{\small time-varying flow trace by finite differences, $308/374/471$ at steps 700, 1500 and 3000} \\
\bottomrule
\end{tabular}
\end{table*}

\subsection{Trajectory Behaviour}

\begin{figure}[!t]
\centering
\includegraphics[width=\columnwidth]{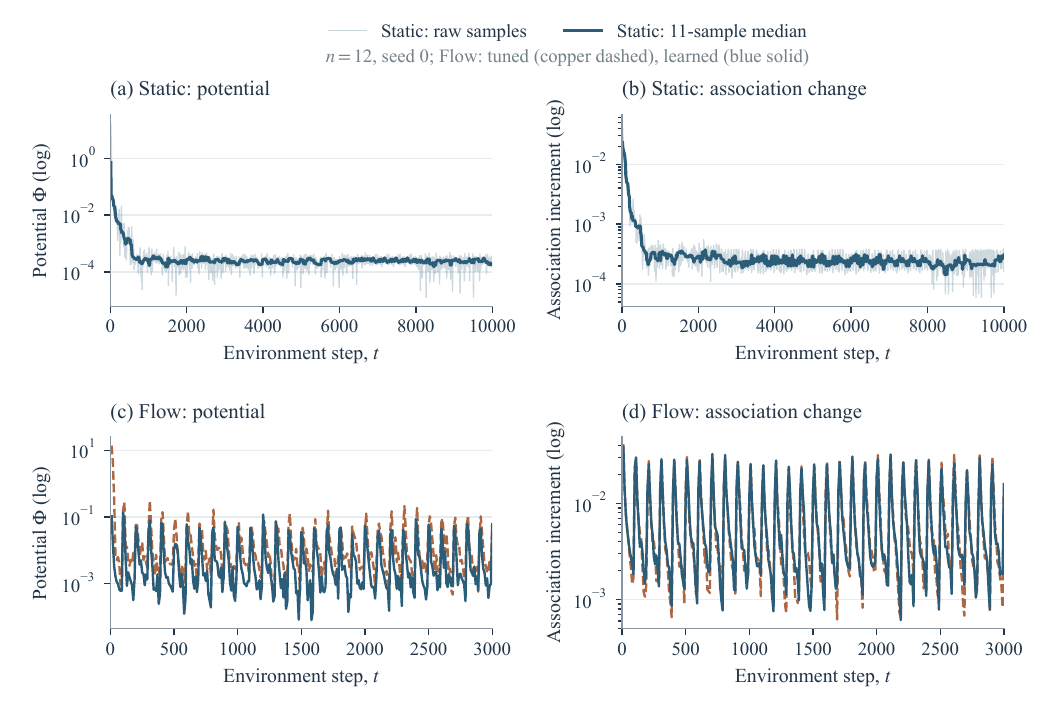}
\caption{Potential and association-update trajectories. On the static instance, left, the potential falls five orders of magnitude on a logarithmic axis, from $17.29$ to $1.31\times10^{-4}$, onto the anchor band. On the time-varying instance, right, the potential tracks inside a bounded band and the association residual on the secondary axis stays at the $10^{-2}$ level.}
\label{fig:traj}
\end{figure}

On the static instance with $n=12$ under the fixed-step executor the potential falls from $17.29$ to the $10^{-4}$ level, as Fig.~\ref{fig:traj} shows. The snapshots at steps 3000, 6000 and 10000 give $\Phi=3.27\times10^{-4}$, $3.52\times10^{-4}$ and $1.31\times10^{-4}$, with load variances $4.53\times10^{-6}$, $4.88\times10^{-6}$ and $1.82\times10^{-6}$ and mean loads steady at $0.7995$, $0.8005$ and $0.7996$ against the anchor of $0.8$, while the association residual decreases monotonically through $3.66\times10^{-4}$, $2.72\times10^{-4}$ and $1.62\times10^{-4}$. The learned executor shows the same diagnostic shape as the fixed-step executor on that instance, since the direction has one source and the speed changes only the time of arrival.

On the time-varying flow instance under the learned executor the potential tracks inside the band $3.6\times10^{-2}$ to $6.1\times10^{-2}$, with the snapshots at steps 700, 1500 and 3000 giving $4.80$, $3.61$ and $6.07\times10^{-2}$, and the association residual stays at the $10^{-2}$ level, with $1.56$, $1.26$ and $1.59\times10^{-2}$. The target moves with the crowd, so the system enters the tracking ball of Proposition~\ref{prop:track} and does not converge to a stationary point. The potential at step 3000 exceeds that of the two earlier snapshots, which corresponds to the crowd roaming into a configuration less favourable to balance, and the system responds inside the band.

Table~\ref{tab:traj} reports the opening trajectory of the time-varying instance step by step. The potential falls two orders of magnitude over the first ninety steps, from $1.05\times10^{-1}$ to $6.18\times10^{-4}$, and the mean load returns from 2.66 to the anchor band at 0.797. At step 100 the potential jumps to $1.35\times10^{-1}$ as a roaming target switch abruptly worsens the balanced configuration, and the system then falls back. Over the whole trace there are 101 rises, the minimum potential is $8.14\times10^{-5}$ and the final value is $6.07\times10^{-2}$. The rises come from target motion, and after each rise the system presses the state back into the band, which is the measured form of the statement of Proposition~\ref{prop:track} that a perturbation enters the tracking ball while the centre of the ball moves with the target.

\begin{table}[!t]
\centering
\caption{Opening trajectory of the time-varying instance with flow, $n=12$ and the learned executor, taken step by step without smoothing.}
\label{tab:traj}
\footnotesize
\setlength{\tabcolsep}{4pt}
\begin{tabular}{@{}rrrrrr@{}}
\toprule
$t$ & $\Phi$ & var$_u$ & $\bar u$ & $a_{\max}$ & dA1 \\
\midrule
10 & $1.05\times10^{-1}$ & $1.45\times10^{-3}$ & 2.661 & 0.976 & $4.01\times10^{-2}$ \\
20 & $1.48\times10^{-2}$ & $2.06\times10^{-4}$ & 1.562 & 0.997 & $1.44\times10^{-2}$ \\
30 & $3.23\times10^{-3}$ & $4.49\times10^{-5}$ & 1.158 & 1.000 & $7.66\times10^{-3}$ \\
40 & $1.66\times10^{-3}$ & $2.30\times10^{-5}$ & 0.961 & 1.000 & $5.35\times10^{-3}$ \\
50 & $1.12\times10^{-3}$ & $1.56\times10^{-5}$ & 0.868 & 1.000 & $4.22\times10^{-3}$ \\
60 & $1.00\times10^{-3}$ & $1.39\times10^{-5}$ & 0.822 & 1.000 & $3.24\times10^{-3}$ \\
70 & $6.75\times10^{-4}$ & $9.38\times10^{-6}$ & 0.797 & 1.000 & $2.78\times10^{-3}$ \\
80 & $6.10\times10^{-4}$ & $8.47\times10^{-6}$ & 0.796 & 1.000 & $2.32\times10^{-3}$ \\
90 & $6.18\times10^{-4}$ & $8.59\times10^{-6}$ & 0.797 & 1.000 & $1.84\times10^{-3}$ \\
100 & $1.35\times10^{-1}$ & $1.88\times10^{-3}$ & 0.909 & 1.000 & $2.36\times10^{-2}$ \\
110 & $7.03\times10^{-2}$ & $9.76\times10^{-4}$ & 0.992 & 1.000 & $3.02\times10^{-2}$ \\
120 & $5.53\times10^{-3}$ & $7.69\times10^{-5}$ & 0.922 & 1.000 & $1.14\times10^{-2}$ \\
\bottomrule
\end{tabular}
\end{table}

\subsection{Scale Axis and Pure-Field Behaviour}

\begin{figure}[!t]
\centering
\includegraphics[width=\columnwidth]{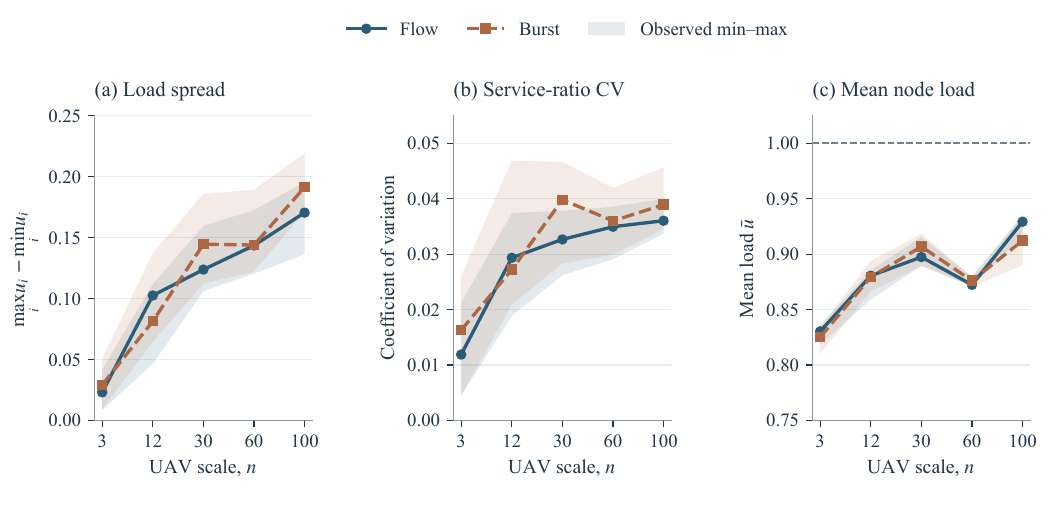}
\caption{Load spread and served-ratio records across scale. The maximum minus minimum load and the mean-field residual stay bounded along the scale axis, with dashed lines marking medians. The mean load at the snapshots of steps 700 and 1500 lies in the band from 0.845 to 0.996.}
\label{fig:scale}
\end{figure}

\begin{figure}[!t]
\centering
\includegraphics[width=\columnwidth]{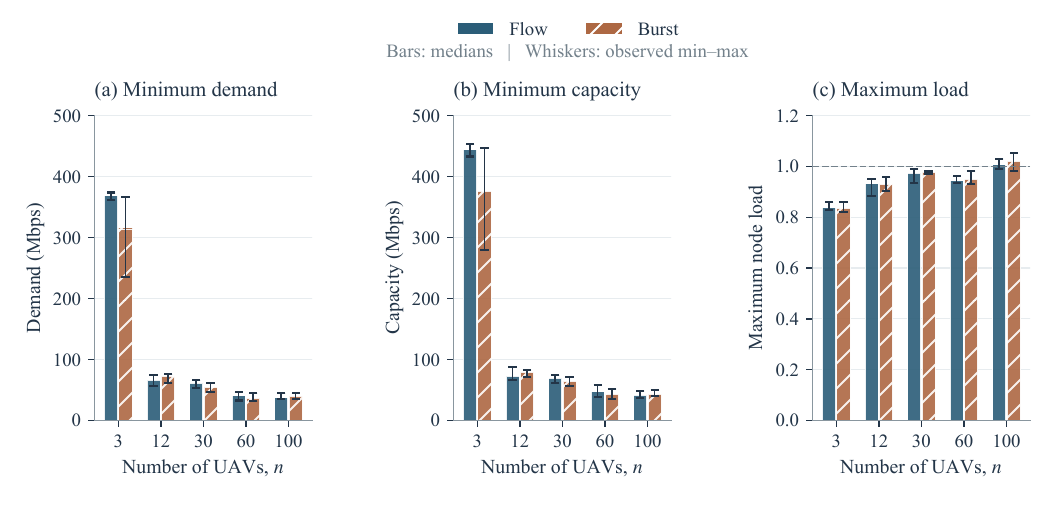}
\caption{Resource and peak-load records, giving the minimum demand, the minimum available rate and the mean load of each scale on the horizontal axis. Resource quantities keep positive lower bounds and the snapshot mean load stays below the overload line.}
\label{fig:resource}
\end{figure}

The sixty states cover three to one hundred UAVs, as Figs.~\ref{fig:scale} and \ref{fig:resource} show. The snapshot mean load lies between $0.845$ and $0.996$, as Table~\ref{tab:spectral} reports, and the minimum demand and the minimum available rate keep positive lower bounds, as Fig.~\ref{fig:resource} shows. Homogeneity, the shape of saturation and the order of the spectral quantities show no trend along the scale axis, so structure is decoupled from scale, which agrees with the complexity analysis of \cite{r69}. Closed-form targets are computed per UAV, broadcast quantities are $O(1)$, and scale only replicates local instances.

The peak number of dead UAVs must be read against the setting. Under the pure-field setting of the audit, with $\eta\equiv1$ and the executor held at a constant gain, the dead count at the snapshot of the burst aggregation transition can reach the whole fleet, with the peak equal to $n$ in every setting of at least twelve UAVs, and later snapshots recover to health, with the mean load back in the band from 0.85 to 1.0. That is the theoretical counterpart of the collapse of the constant-speed scheme in Table~VII of \cite{r69}. When the direction is correct but the speed cannot adapt, the demand shock of an instant of aggregation pushes part of the fleet into deep overload, and a maximal-speed catch-up amplifies oscillation in the ill-conditioned region of Proposition~\ref{prop:beta}. The structural checks of homogeneity, energy decomposition and spectral quantities hold on those very states, so structural properties and execution behaviour live on separate layers. The former belong to the closed-form direction and always hold, while the latter depends on the speed law. Under the learned executor the same scenarios of the main experiment of \cite{r69} show almost no dead UAVs, with twenty-eight of thirty settings at zero in Table~VI of \cite{r69}, and the two sets of numbers together give the complete picture. The pure field supplies the direction and the structure, and the learned gain pulls execution behaviour from a fleet-wide peak back into the zero-death band.

\subsection{What the Audit Establishes and What It Leaves Open}

The audit is a measurement, not a proof, and the boundary between the two is worth stating before the numbers. Three of the four groups are checks of structural identities that the model makes exact: homogeneity of the load in bandwidth, the energy identity of the conservative gain map, and the anchoring of the position target at the demand-weighted centroid. Each of these is verified to numerical precision on measured states, so a discrepancy in any of them would falsify the structure and not merely its constants. The fourth group is different in kind: the spectral radius is a property of the linearisation at a point, and the audit reports its value across operating points, which is a statement about how far the local behaviour may differ from the global claim.

The audit therefore supports two conclusions and neither more. The structural claims hold on the assembled pipeline, so they are not artefacts of a simplified model. The spectral claim is a statement about the one-step association map measured on states of the audit, so it can differ from the spectral radius of the composite that runs in the main experiment. The difference is not a defect and it is not a contradiction: the certificate is used to show that a single global contraction rate does not exist, and that conclusion only requires the amplification to exceed one somewhere, which the audit shows over a wide range of operating points.

Two consequences follow for how the rest of the paper should be read. The convergence statements of Section~IV do not depend on the audit, since they are proved from the one-dimensional map. The empirical statements of Section~V do not depend on the spectral radius either, since they compare measured indicators. The audit sits between the two and explains why the layered certificate has the form it has.

\section{Discussion}

The coverage of the guarantees is stated layer by layer.

The structural layer is unconditional. The homogeneous decomposition, the anchored closed form and centroid optimality are algebraic identities or strict-convexity conclusions that hold for every channel distribution, every operating point and every parameter value, at every instant of the system operation, and the Lipschitz and convex-domain premises of a general resource allocation analysis are superseded here by identities. The audit checks sixty states with no exception, which is the empirical counterpart of that universality.

The execution layer covers the entire gain domain. Theorem~\ref{thm:1d} delivers a complete four-part characterisation over the whole space of positive gains, arbitrary initial values and bounded per-step perturbations, namely geometric convergence, the critical two-state oscillation, a bounded family of two-state cycles and global boundedness, with no premise on gain, initial value or perturbation. Its guarantee domain is strictly larger than the applicability domain of spectral compression, since the spectral certificate fails on all sixty measured states at $10^{1}$ to $10^{3}$ while the gain-truncation mechanism stays valid over that same region. The rate is characterised by the distinction between the geometric phase and the cycle phase together with an explicit step ceiling from Corollary~\ref{cor:steps}, and every part of that statement holds over the whole gain domain.

The multidimensional composition follows the division into safe per dimension and coupled through the perturbation. The target motion caused by environment coupling enters $\delta_{\max}$, bounded term by term in Table~\ref{tab:disturb}. Any gain keeps the action bounded, and bounded target increments together with bounded disturbances give the explicit ultimate bound of Proposition~\ref{prop:track}; that bound, and not the size of a single step, is what the operating point of a moving scene is measured against.

The audit and the main experiment complement each other across layers. The audit, run in the pure field with a fixed channel class, characterises structure through homogeneity, energy decomposition, spectral quantities and pure-field behaviour. The main experiment of \cite{r69}, run with the learned executor under the primary channel setting, characterises post-training behaviour, with almost no dead UAVs across thirty settings and survival in every high-pressure setting. The two layers together form a complete chain of evidence. The structural layer explains why the behavioural layer holds, through a closed-form direction and a safe gain, and the behavioural layer confirms the design value of the structural layer, since the fleet-wide peak of the pure field is pulled back into the zero-death band by the learned gain.

\subsection{Complexity and Deployment}

The cost of the analysed system is local in three senses. The decider evaluates a closed form per UAV over the users it serves, so its work is proportional to the number of served users and to the size of the state, and no matrix over the fleet is formed. The executor evaluates a network whose input and output dimensions are fixed by the action dimension, so its cost does not grow with $n$. The only fleet-wide quantities are two scalars, the sum of the demand scales and the mean price, both refreshed on a slow clock, which means the coupling across the fleet costs two broadcasts per refresh, with no global solve.

This is the sense in which the layered structure delivers scalability that a centralised optimiser cannot. A block-coordinate method re-solves a bandwidth program and a trajectory program over the whole fleet on a duty cycle, and its cost grows with the number of UAVs and users together. The analysis here shows that the object that such a program would compute, the equal-load bandwidth target, has a closed form once the load is written in the right view, so the expensive part can be replaced by an identity. What remains centralised is the evaluation of a fleet-wide potential for reporting, which the deployment can compute offline.

The guarantees translate into deployment statements without further assumptions. Because the gain is bounded and the action is clipped, the executed action never leaves the physical envelope, so the safety of the platform is a property of the action space and not of the learned weights. Because the direction has no parameters, a deployment can run the decider with a constant gain before any training and improve from there, and the early-stop band of the training procedure is defined relative to that starting point. Because the tracking bound is explicit in the disturbance magnitudes, an operator can read the expected ultimate error from the refresh period of the position information and the scale of the demand motion, which are the two quantities that a deployment knows.

A final deployment point concerns the audit. The structural properties are verified on measured states, not only on the model, and the verification is a small script over the states of a run. The same script runs on the states of a new deployment, which turns the audit into a commissioning check: homogeneity, the energy identity, the potential and the position target can be confirmed on the deployed fleet before the gains are tuned.

\section{Conclusion}

This paper proves the core mathematical properties of the companion system \cite{r69}. Under demand-proportional bandwidth allocation the load admits a first-order homogeneous decomposition in bandwidth together with an anchored closed form, including feasibility bounds from the budget and the backhaul. The logarithmic knee view preserves order and characterises the zero set of the potential. The position potential has a unique minimiser at the weighted centroid. The potential and the Jain fairness index share their zero set, are monotone in the same dispersion on the linear band of the view, and together define the equal-load optimum. The association block possesses a fixed point. The one-dimensional truncated gain iteration is characterised completely, through geometric convergence, an explicit step ceiling, a two-state bounded cycle and global boundedness. A perturbed tracking ball and the step ceiling $\eta^\ast\propto\beta^{-1/2}$ complete the set.

The numerical audit verifies homogeneity in sixty of sixty states. The Euclidean spectral contraction condition fails on all sixty measured states while the system converges everywhere and survives in every setting, which shows that the combination of a hard gain ceiling and truncation delivers guarantees over the entire operating region in which the current spectral methods fail. That is the substantive strength of this framework over a spectral route. The step must adapt to the ill-conditioning gain, which is precisely the mathematical ground of an executor that learns only the speed while the direction stays closed form.

Two extensions follow naturally. A joint theorem for the convergence of the outer alternation between the association block and the executor would close the composite iteration, and the static trace, where the potential decreases monotonically by five orders of magnitude from $17.29$ to $1.31\times10^{-4}$, together with the critical spectrum near one already supports that direction numerically. A quantitative relation between the amplitude of the two-state family and the return-optimal gain would lift the present safety characterisation into an optimality characterisation.

\appendices

\section{Numerical Method}

\subsection{Construction of the Audit Probes}

The homogeneity probe D1 applies a relative bandwidth perturbation in both directions to every UAV at a snapshot, recomputes the whole chain under the frozen association condition, and compares the drift of $c_i=u_ib_i$, reporting the residual as a relative quantity. The spectral probe D2 estimates the largest eigenvalue in the Euclidean norm of the Jacobian of the one-step association map, taken at frozen position, frozen price and frozen demand, with the bandwidth shares, the transmit power, the interference denominator and the available rate recomputed from the perturbed association at every matrix-vector product, so that the derivative carries the diagonal relaxation of the damped softmax together with the response of the whole chain. The matrix-vector product is formed by forward differences at a snapshot and the resulting operator is iterated with ARPACK for the largest eigenvalue in the Euclidean norm. The difference step is swept over four values from $10^{-4}$ to $10^{-10}$, and the spread of the four estimates at one point reports the stability of the estimate. That spread stays within $0.3\%$ inside one snapshot, going from $1.0246$ to $1.0268$ at step 3000, from $0.9801$ to $0.9792$ at step 6000 and from $0.9953$ to $0.9950$ at step 10000, while the differences between snapshots, $1.027$, $0.979$ and $0.995$, come from the different snapshot points and not from numerical noise. The convergence diagnostic records the potential, the variance, the mean load, the association residual, the largest association weight and the successive difference of the KL quantity at every step for a designated instance, with the full set of spectral estimates added at snapshots.

\subsection{Scale and Numerical Setting of the Audit}

The audit covers five scales, two scenario classes and three seeds, giving thirty instances and sixty snapshot states, one at step 700 and one at 1500. Three instances, namely a static one under the fixed-step executor, a time-varying one under the fixed-step executor and a time-varying one under the learned executor, are additionally recorded step by step for the potential, the variance, the mean load and the spectrum, which shows how the shape of convergence varies with scale and pressure. The audit shares the environment code, the link computation and the demand model with the main experiment, and the three differences listed below are what keep the structural probes independent of the trained weights. The results of the audit and those of the companion system \cite{r69} are therefore read each against its own setting, and no number is carried between the two as a controlled ablation.

\subsection{Setting of the Audit}

The audit differs from the main experiment of \cite{r69} in three respects, each with a purpose. It executes the pure field with $\eta\equiv1$, so that the structural checks are completely independent of the trained weights and hold for any network weights. It uses the UMi channel class with a fixed 6 dB shadow, which removes the randomness between episodes and makes spectral estimates comparable across snapshots. It runs without the physical engine, so every probe acts on the algebraic chain alone and no airframe response enters a measurement. The main experiment performs every performance evaluation under the learned executor, the primary channel setting and the full physical engine. Two labels are kept apart throughout: pure field means $\eta\equiv1$ in the audit, and constant speed means the always-maximal gain baseline of the main experiment. They are different executors and their numbers are not interchangeable. Homogeneity is an algebraic identity and Lemma~\ref{lem:homogeneous} is valid for any channel class, so the measured UMi result and the proof confirm each other.

\section{Derivation Details of the Ill-Conditioning Gain and the Step Bound}

This appendix expands the explicit bounds of the individual links in the proof of Proposition~\ref{prop:beta}, all of which follow by differentiating the model.

\subsection{Sensitivity of the Spectral Sum to Spectral Efficiency}

From (\ref{eq:Si}) one has $S_i=\sum_k(x_{ik}^2d_k/\rho_i)r_{ik}$, and differentiating term by term while using $x_{ik}\le1$ gives
\begin{equation}
\frac{\partial S_i}{\partial r_{ik}}=\frac{x_{ik}^2d_k}{\rho_i}\ \le\ \frac{x_{ik}d_k}{\rho_i},\qquad
\sum_k\frac{x_{ik}d_k}{\rho_i}=1.
\end{equation}
As a function of $r$, the spectral sum thus has a gradient whose one-norm stays at most one, so it does not amplify spectral efficiency. The demand scale $c_i=\rho_i/S_i$ is relatively sensitive to $S_i$, with sensitivity $1/S_i$, since a ten per cent drop in $S_i$ raises $c_i$ by about eleven per cent, and $1/S_i^2$ therefore enters the ill-conditioning gain.

\subsection{Sensitivity of the Association Logit to Demand}

Differentiating the rate-margin term of $z_{ik}$ with respect to $\rho_i$ at fixed $R$ gives
\begin{equation}
\frac{\partial}{\partial\rho_i}\frac{R_i^{\mathrm{avail}}}{\rho_i+\varepsilon}=-\frac{R_i^{\mathrm{avail}}}{(\rho_i+\varepsilon)^2},
\end{equation}
so the gradient flattens as demand grows, through the squared denominator, while small demand is regularised by $\varepsilon$ and keeps the derivative bounded by $R/(\rho+\varepsilon)^2$. Multiplication by $\kappa$ yields the second entry of $\beta$.

\subsection{Gain From Bandwidth to Load}

Differentiating $u_i=c_i/b_i$ with respect to $b_i$ in the access-limited region gives
\begin{equation}
\frac{\partial u_i}{\partial b_i}=-\frac{c_i}{b_i^2}=-\frac{u_i}{b_i}=-\frac{u_i^2}{c_i}.
\end{equation}
The two ways of writing one derivative expose the ill conditioning at both ends. The form $u_i/b_i$ grows when bandwidth is small, and the form $u_i^2/c_i=\rho_i/(S_ib_i^2)$ grows with the carried demand and as the spectral sum falls.

\subsection{Composition}

Multiplying the links of the preceding subsections with the well-conditioned links, namely the view with $m'\le1$, the softmax with its $\tfrac14$-Lipschitz property \cite{r63} and the clip with its $1$-Lipschitz property, and merging the common factor into $\sqrt\beta$ with $\beta$ as defined in Proposition~\ref{prop:beta}, gives the upper bound $L_F\le\eta C_\ast\sqrt\beta$ on the closed-loop Lipschitz constant, where $C_\ast$ contains only well-conditioned constants. The contraction condition $L_F<1$ yields the step ceiling $\eta^\ast=C_\ast^{-1}\beta^{-1/2}$. A tenfold increase of $\beta$ contracts the admissible step to about one third of its former value, a factor of $\sqrt{10}\approx3.2$. The bound is a sufficient condition: below $\eta^\ast$ contraction follows, and beyond it no conclusion is drawn. The mathematical consequence is that one fixed step cannot satisfy the condition over a wide range of $\beta$, since a step small enough at high $\beta$ leaves the iteration slow at low $\beta$. The experimental evidence for gain adaptivity is separate and does not rest on this bound: the same-direction fixed-step executor collapses in the deep-overload grids of \cite{r69}, where the spectral sum falls and the carried demand rises at once, while the learned gain recovers performance on the same grids.

\section{Standing Assumptions}
\label{app:assumptions}

The propositions name their hypotheses, and this section gathers the model assumptions that the analysis of the companion system fixes. They are listed once so that each statement can be read against them.

\begin{enumerate}
\item[\textbf{A1}] \emph{Demand-proportional split.} Inside a UAV the bandwidth of its associated users is divided in proportion to $x_{ik}d_k$, so the share of user $k$ is $b_i x_{ik}d_k/\sum_j x_{jk}d_j$. The access rate of a UAV is the sum over its users of the share times the spectral efficiency of the link.
\item[\textbf{A2}] \emph{Access-limited region.} The closed forms of Section~III are derived where the access rate is below the backhaul capacity. Outside that region the backhaul binds, the per-user rates are scaled by a common factor, and the load is that of the capped rate. The cap is a structural part of the model and is not assumed to bind in the main experiment.
\item[\textbf{A3}] \emph{Load view.} The view is the logarithmic knee map, the identity for $u\le3$ and $3+3\ln(1+(u-3)/3)$ above it. It is continuous, strictly increasing and concave above the knee, so order is preserved and no load is mapped to a common ceiling.
\item[\textbf{A4}] \emph{Association damping.} The association follows a damped softmax with damping time $\tau=5$ and logits $\alpha r_{ik}+\kappa R_i/(\rho_i+\varepsilon)+\pi_i$ with $\alpha=1$, $\kappa=0.5$ and $\varepsilon=0.1$. The logits depend on the current rates, which depend on the association, so the block is a map on the product of simplices and not an optimisation over $X$.
\item[\textbf{A5}] \emph{Information latency.} The direction decider sees user positions that refresh every $0.5$ s with ten metres of Gaussian noise and holds them between refreshes, while the environment advances on true positions. The position target is therefore a perturbed centroid, and the perturbation enters the tracking bound of Section~IV: the deterministic recursion takes its size at the ten-metre scale, and the Gaussian tail beyond that scale is what the high-probability form of the bound accounts for.
\item[\textbf{A6}] \emph{Gain and clip.} The learned gain of each dimension lies in $[0.497,2.014]$, the action is clipped at the physical single-step bound, and the bandwidth and price dimensions carry an additional relaxation of $0.2$. The clip is part of the map that Theorem~1 analyses, so boundedness is a property of the implemented system and not of an idealised one.
\end{enumerate}

Assumption~A1 is what makes the load first-order homogeneous in bandwidth, A2 fixes where the closed form applies, and A4 is what rules out a global gradient argument: the association is a fixed-point map, so the composite dynamics are analysed in the neighbourhood of its fixed point, with no descent on the potential. Assumptions~A5 and A6 supply the two bounded quantities that the tracking recursion consumes.

\section{Worked Orbits and Boundary Arithmetic}
\label{app:orbits}

The statements of Theorem~1 are one-dimensional, and the arithmetic behind each branch can be checked by hand. This appendix carries the arithmetic at unit single-step bound, $h=1$, and at a target placed at the origin, so the state and the error coincide.

\subsection{Below the Critical Gain}

At $\gamma=1.5$ the map is the linear one wherever $|e|\le2/3$, and the saturation is reached only at the start. From $e_0=0.9$ the first step saturates and lands at $e_1=0.9-1=-0.1$, which is inside the linear region, and the following steps contract by the factor $|1-\gamma|=0.5$: the sequence is $-0.1$, $0.05$, $-0.025$, $0.0125$, and so on. The sign flips at every step because the gain overshoots, and the magnitude halves, which is the geometry of the geometric branch of the theorem.

\subsection{At the Critical Gain}

At $\gamma=2$ the linear region is $|e|\le1/2$, and inside it the map is $e\mapsto-e$, a pure reflection. From $e_0=0.9$ the first step saturates to $e_1=-0.1$, the gain applied to $-0.1$ has magnitude $0.2$, which stays below the bound, so the second step is linear and gives $e_2=+0.1$, and from there the orbit alternates between $0.1$ and $-0.1$. The pair $\{0.1,-0.1\}$ is the two-state cycle of the critical case, its magnitude is set by the initial value, not by the gain, and it sits below the bound $h/2$. An initial value inside the linear region never leaves it, and an initial value with $|e_0|\le h$ is mapped into that region in one step, while a larger initial value descends by exactly $h$ per saturated step until it lands inside.

\subsection{Above the Critical Gain}

At $\gamma=2.5$ the cycle family covers $p\in[0.4,0.6]$. From $e_0=0.9$ the orbit runs
\begin{equation*}
\begin{aligned}
0.9&\ \to\ -0.1\ \to\ 0.15\ \to\ -0.225\ \to\ 0.3375\\
&\ \to\ -0.50625\ \to\ 0.49375\ \to\ -0.50625,
\end{aligned}
\end{equation*}
where the fifth and sixth steps both saturate, and the orbit then repeats the pair $\{0.49375,-0.50625\}$. The first element of the pair lies in $[0.4,0.6]$, as the family requires, and the magnitude of the pair is $0.50625$, below the ceiling $(1-1/\gamma)h=0.6$. Every state of the orbit after the first step stays inside that ceiling, and the largest magnitude of the whole orbit is the negative element of the cycle.

At $\gamma=3$ the family covers $p\in[1/3,2/3]$, and from the same initial value the orbit reaches the pair $\{0.6,-0.4\}$, with the positive element at the upper end of the family and a magnitude of $0.6$, below $2/3$. In both cases the two elements of the cycle differ by exactly $h$ and the positive element lies inside the family interval, while neither element sits at an end of that interval: the interval is the set of admissible positive elements, not a range that the orbit sweeps, and which element of it the orbit settles on is fixed by the initial value.

\begin{table}[!t]
\centering
\caption{Orbits at unit single-step bound from $e_0=0.9$.}
\label{tab:orbits}
\footnotesize
\begin{tabular}{@{}cccl@{}}
\toprule
$\gamma$ & Regime & Limit object & Magnitude bound \\
\midrule
$1.5$ & geometric & target & $|1-\gamma|^t|e_0|$ \\
$2$ & critical & $\{0.1,-0.1\}$ & $h/2$ \\
$2.5$ & supercritical & $\{0.494,-0.506\}$ & $(1-1/\gamma)h=0.6$ \\
$3$ & supercritical & $\{0.6,-0.4\}$ & $(1-1/\gamma)h=0.667$ \\
\bottomrule
\end{tabular}
\end{table}

\subsection{Boundary of the Closed Form}

The feasibility level of Proposition~3 has two terms, and the branch that binds depends on the instance. For a fleet in which the demand scales are equal, $c_i=c$, the budget term reads $nc/B$ and the backhaul term reads $c^{\mathrm B}$, so the budget binds when $nc/B\ge c^{\mathrm B}$, which is the regime of the main experiment at large scale. For a fleet with one heavily loaded UAV, the maximum over $i$ of $c_i^{\mathrm B}$ carries the binding role and the closed form saturates at the backhaul. The two branches meet at equality, and the target is not unique there, since any level above $\mu_{\mathrm{low}}$ satisfies both constraints. The proposition therefore claims a minimal feasible level and not a unique one, and the algorithm uses the anchor $u^{\ast}=0.8$ whenever it exceeds the level.

\section{Proof Details}
\label{app:proofs}

The steps below fill in the algebra behind the propositions of Section~III and Section~IV. Throughout, $h$ is the single-step bound of the implementation, $\eta$ the gain, $\gamma=\eta\beta$, and $T$ the map of \eqref{eq:1d}, namely $T(e)=e-\mathrm{clip}(\gamma e,\pm h)$.

\subsection{Shape of the One-Dimensional Map}

The map is a composition of the linear gain $e\mapsto\gamma e$, the clamp onto $[-h,h]$ and the subtraction, and each of the three is monotone nondecreasing, so $T$ is monotone nondecreasing. Its three pieces follow by asking whether the clamp is active. When $|\gamma e|\le h$, that is $|e|\le h/\gamma$, the clamp is the identity and $T(e)=(1-\gamma)e$. When $e>h/\gamma$ the clamp returns $h$ and $T(e)=e-h$, whose slope is one. When $e<-h/\gamma$ the clamp returns $-h$ and $T(e)=e+h$, again of slope one. The map is odd, since $T(-e)=-T(e)$ for every $e$.

The three pieces meet continuously at the junctions $e=\pm h/\gamma$, where $e-h$ and $(1-\gamma)e$ agree, as do $e+h$ and $(1-\gamma)e$: at $e=h/\gamma$ both give $h/\gamma-h$, and at $e=-h/\gamma$ both give $h-h/\gamma$. On the slab $[-h,h]$ the image is therefore the interval $[-(1-1/\gamma)h,(1-1/\gamma)h]$, which lies inside the slab. The map has one hump on the right, one valley in the middle and one rise on the left, and the extreme values of the image at the two clamped ends are $T(h)=0$ and $T(-h)=0$.

The zeros of $T$ are the solutions of $T(e)=0$: the linear piece gives $e=0$, the right piece gives $e=h$ and the left piece gives $e=-h$. The value $e=h$ is admissible because at $e=h>h/\gamma$ the clamp saturates and a step of exactly $h$ lands on the target, which is the one-step arrival used in the first claim of Theorem~\ref{thm:1d}. There is no other preimage of the target inside the slab.

\subsection{The Saturated Step Never Overshoots the Target Twice}

From $h/\gamma<e\le h$ the image is $e-h\le0$, whose magnitude is $h-e<h-h/\gamma$. Since $\gamma<2$ gives $h-h/\gamma<h/\gamma$, that magnitude lies in the linear range and the following step is linear, which is the mechanism that keeps the orbit from alternating in the clipped region when the gain is below two. From $e>h$ the image is $e-h>0$ and the magnitude decreases by exactly $h$ per step until the trajectory enters $[0,h]$, at which point the previous case applies. The symmetric sentence holds on the negative side. This is the content of the first claim of Theorem~\ref{thm:1d}, and it makes no use of any monotonicity of the potential, only of the arithmetic of the clamp.

\subsection{Verification of the Two-State Family}

Let $\gamma>2$ and let $p\in[h/\gamma,(1-1/\gamma)h]$. The upper end of the interval is equivalent to $\gamma(h-p)\ge h$, and hence to $h-p\ge h/\gamma$. For the point $p$ the clamp saturates, so $T(p)=p-h$. For the partner $q=p-h$ one has $q<0$ and $|q|=h-p\ge h/\gamma$, so the left clamp saturates and $T(q)=q+h=p$. The pair $\{p,p-h\}$ is a two-cycle. Its amplitude is $\max\{p,h-p\}$, and since $p\le(1-1/\gamma)h$ and $h-p\le h-h/\gamma$, the amplitude is at most $(1-1/\gamma)h$, which is below $h$ for every $\gamma>2$ and approaches $h/2$ as $\gamma$ approaches two from above.

At $\gamma=2$ the interval shrinks to the single point $p=h/2$ and the pair is $\{h/2,-h/2\}$. A linear step at $\gamma=2$ is $e\mapsto-e$, so once the orbit lands in $|e|\le h/2$ it alternates with a magnitude fixed by its initial value, and if it lands exactly on the target it stops. Both endings are admitted by the second claim of Theorem~\ref{thm:1d}.

\subsection{Uniqueness of the Two-State Family}

A two-cycle of $T$ that is not the pair of the family would have to be produced by two different branches, because the linear branch is injective on its range and the two clipped branches are injective. A cycle that uses the linear branch twice would be a fixed point of a contraction, hence the fixed point itself, and a cycle that uses one linear and one clipped branch would require $(1-\gamma)e$ to land exactly on the junction, which happens only at the critical gain. The only remaining possibility is the pair of clipped branches, which is the family just computed.

For the convergence claim, a bounded orbit that enters $|e|\le h/\gamma$ is multiplied by $\gamma-1>1$ on every pass through the linear piece, so it leaves that range after finitely many steps and saturates again; hence it cannot converge to the target. A bounded orbit that never enters the linear piece lives in the two clipped regions, where the composition of the two branches is the identity, so it is one of the family pairs. This is the third claim.

\subsection{Global Boundedness}

Let $V_t=|e_t|$. From $V>h$ the saturated step gives $V_{t+1}=V-h<V$. From $h/\gamma<V\le h$ the saturated step gives $V_{t+1}=h-V<V$, and since $V>h/\gamma$ we have $V_{t+1}<h-h/\gamma<h$. From $V\le h/\gamma$ the linear step gives $V_{t+1}=\gamma' V$ with $\gamma'=|\gamma-1|$, which satisfies $\gamma'V\le\gamma'h/\gamma<h$ because $|\gamma-1|<\gamma$ for every positive $\gamma$. Induction on $t$ gives $V_t\le\max\{V_0,h\}$ for all $t$, and every step lands in the slab, which is therefore positively invariant. This is the fourth claim.

\subsection{Homogeneity Algebra}

In the access-limited region the load of UAV $i$ is $u_i=\rho_i/(b_iS_i)$, and the budget share of user $k$ at UAV $i$ is $x_{ik}/\sum_jx_{jk}$, so scaling every bandwidth by $\lambda$ scales every $R_i^{\mathrm A}$ by $\lambda$ and leaves every assignment $x_{ik}$ unchanged. The exponent of the backhaul hop does not enter, since the region assumes $b_iS_i\le C_i^{\mathrm B}$. Hence $u_i(\lambda b)=\lambda^{-1}u_i(b)$ throughout the region, and the target $b_i^\ast=c_i/u^\ast$ satisfies $\lambda b_i^\ast(\lambda^{-1}u^\ast)=b_i^\ast(u^\ast)$, which is homogeneity in the sense of Section~III. The anchored potential is invariant because it is a function of the normalised loads and of centroid distances, and neither changes under the rescaling.

\subsection{Centroid Optimality}

With the position potential of Section~II-B, the derivative with respect to the position of UAV $i$ is
\begin{equation*}
\frac{\partial\Phi_P}{\partial\xi_i}=2\sum_k w_k(\xi_i-y_k),\qquad w_k=\frac{c_ix_{ik}}{\rho_i\,|\xi_i-y_k|^{2}},
\end{equation*}
so the stationary point is the weighted mean $\xi_i^\ast=\sum_kw_ky_k/\sum_kw_k$. The weight is the share of the user demand in the load of the UAV discounted by the inverse square distance, which makes the stationary point a demand-weighted centroid of the served users. In the balanced state every distance is equal, the weights reduce to $c_ix_{ik}/\rho_i$, and the stationary point becomes the demand-weighted centroid of Proposition~3. The Hessian is positive definite whenever the weights are positive, so the stationary point is the minimum, and because the objective is separable across coordinates and the constraint set is a box, the projected iterate is the clipped move towards that point.

\subsection{Jain Index in the Continuous Limit}

For $n$ loads with mean $\bar u$ the index is $(\sum_iu_i)^{2}/(n\sum_iu_i^{2})$, which equals one when $u_i=u_j$ for all pairs and falls towards $1/n$ when one load carries the whole sum. Substituting $u_i=\rho_i/b_i$ gives the expression in Proposition~5, and the identification of the maximum with the equal-load configuration follows from Cauchy--Schwarz. The potential shares the maximiser because its unique zero is the equal-load configuration of the same instance.

\subsection{The Damping Estimate}

The association recursion is $a_{t+1}=(1-\tau^{-1})a_t+\tau^{-1}x_t$ with $\tau=5$. Subtracting a fixed point $a_\ast$ gives $e_{t+1}=(1-\tau^{-1})e_t+\tau^{-1}(x_t-x_\ast)$, and taking norms yields
\begin{equation*}
\|e_{t+1}\|\le(1-\tau^{-1})\|e_t\|+\tau^{-1}\|x_t-x_\ast\|,
\end{equation*}
so the block carries the geometric factor $1-\tau^{-1}=0.8$ on its own and a forcing term of size $\tau^{-1}$ times the softmax motion. The factor 0.8 is why the association block is damped even though it is not the source of the amplification reported in Section~V-B.

\section{The Audit Protocol}
\label{app:protocol}

The audit takes the states of a run and asks of each state whether the object that the analysis names is actually present in it. This appendix records how a state is captured, what each check computes, and which tolerances the report uses, so that the same procedure can be replayed on another run.

\subsection{State Capture}

A state is a snapshot of the assembled pipeline at a control step: the served assignments $x_{ik}$, the bandwidths $b_i$, the prices $\pi_i$, the power multipliers $\tilde p_i$, the positions $\xi_i$, the user demands $d_k$ and the channel gains $G_{ij}$ with their fitted exponents. The capture does not interrupt the run: the environment writes the state to a record and continues, and the fields of the record are the same fields that the main experiment stores, so a state from the audit and a state from the main experiment are drawn from the same distribution of records.

Thirty instances give sixty states, two per instance, taken at two points of the trajectory. The early point is inside the settling phase and the late point is after the association block has settled, which is deliberate: the checks are meant to hold while the state is moving as well as when it has stopped.

\subsection{The Four Checks}

\emph{Homogeneity.} The load of each UAV is recomputed from the same assignment under a rescaling of the bandwidth vector by a factor chosen so that the total stays on the budget. The quantity $c_i=u_ib_i$ is compared before and after, and the check passes when the relative deviation is at machine precision. The reported drift over sixty states is $2.7\times10^{-16}$, which is the level of the double-precision round-off of the product.

\emph{Energy identity.} The conservative gain map of the executor carries an energy that is quadratic in the residual of the idempotent map, and the identity states that the decrease across one update equals the dissipated part. The check recomputes both sides from the state and compares them, which verifies that the trunk of the network has not introduced an asymmetric term.

\emph{Potential.} The anchored potential is recomputed from the normalised loads and compared with the value that the environment reports, which verifies that the two implementations agree on the normalisation and on the anchor.

\emph{Position target.} The demand-weighted centroid of the served users is recomputed and compared with the position target that the decider produced, which verifies the closed form of Proposition~3 on measured assignments.

\subsection{The Spectral Measurement}

The association map is measured by a matrix-free power iteration on the linearisation. A small perturbation is applied to the assignment vector of one UAV, the map is evaluated once at the perturbed state, and the difference quotient gives the action of the Jacobian on that direction. The leading eigenvalue is extracted from the resulting operator with a Lanczos-type solver at a tolerance of $10^{-3}$ and a cap of three hundred iterations.\footnote{The implementation uses the sparse linear operator interface of the standard scientific stack, with the perturbation size chosen to stay above the round-off floor of the difference quotient.}

Two properties of the measurement matter for reading its value. The first is that the derivative is taken along the assembled chain, so the number includes the bandwidth share, the power multiplier, the interference denominator and the available rate of every UAV, and not the association block alone. The second is that the perturbation is applied to one UAV at a time, so the number is a local amplification of that UAV's assignment and not a fleet-wide multiplier.

The reported range over the operating points is from fifty to three thousand one hundred, with the smallest values at the lightly loaded ends of the pressure range and the largest where the carried demand sits far from the mid range. The range is what the certificate needs: one value above one somewhere already rules out a single global contraction factor, and the measurement shows that the excess is not marginal.

\subsection{What the Audit Was Not}

The audit does not run the physical engine, does not use the trained executor and does not use the shadowing field of the main experiment. Its purpose is to isolate the algebraic structure, so the gain is held at one and the executor is replaced by its identity, which is the pure-field setting named in Section~V. The two settings therefore differ in three respects and cannot be substituted for one another. The main experiment measures behaviour, the audit measures structure, and the two are reported side by side as complementary evidence rather than as a controlled ablation of one another.

\vspace{6pt}

\noindent\textbf{Open source.} The implementation and the full experimental records are available at \url{https://github.com/jhzAliyy/space4-run}.

\end{document}